\documentclass[12pt]{article}

\usepackage{graphicx}

\scrollmode
\usepackage{amsfonts}
\usepackage{amsmath}
\usepackage{amssymb}
\begin{document}
\title{\bf The Brownian Web: Characterization and Convergence}
\author{L.~R.~G.~Fontes\thanks{Partially supported by CNPq grants
  300576/92-7
and 662177/96-7 (PRONEX) and \quad\mbox{ }
\hspace{.3cm} FAPESP grant 99/11962-9} 
\and
M.~Isopi\thanks{Partially supported by the CNPq-CNR
agreement} \and C.~M.~Newman\thanks{Partially supported by NSF grant DMS-0104278}
\and K.~Ravishankar\thanks{Partially supported by NSF grant DMS-9803267}}
\date{04/08/03}
\maketitle

\newtheorem{defin}{Definition}[section]
\newtheorem{Prop}{Proposition}
\newtheorem{teo}{Theorem}[section]
\newtheorem{ml}{Main Lemma}
\newtheorem{con}{Conjecture}
\newtheorem{cond}{Condition}
\newtheorem{prop}[teo]{Proposition}
\newtheorem{lem}{Lemma}[section]
\newtheorem{rmk}[teo]{Remark}
\newtheorem{cor}{Corollary}[section]
\renewcommand{\theequation}{\thesection .\arabic{equation}}

\newcommand{\beq}{\begin{equation}}
\newcommand{\eeq}{\end{equation}}
\newcommand{\beqn}{\begin{eqnarray}}
\newcommand{\beqnn}{\begin{eqnarray*}}
\newcommand{\eeqn}{\end{eqnarray}}
\newcommand{\eeqnn}{\end{eqnarray*}}
\newcommand{\bprop}{\begin{prop}}
\newcommand{\eprop}{\end{prop}}
\newcommand{\bteo}{\begin{teo}}
\newcommand{\bcor}{\begin{cor}}
\newcommand{\ecor}{\end{cor}}
\newcommand{\bcon}{\begin{con}}
\newcommand{\econ}{\end{con}}
\newcommand{\bcond}{\begin{cond}}
\newcommand{\econd}{\end{cond}}
\newcommand{\eteo}{\end{teo}}
\newcommand{\brm}{\begin{rmk}}
\newcommand{\erm}{\end{rmk}}
\newcommand{\blem}{\begin{lem}}
\newcommand{\elem}{\end{lem}}
\newcommand{\ben}{\begin{enumerate}}
\newcommand{\een}{\end{enumerate}}
\newcommand{\bei}{\begin{itemize}}
\newcommand{\eei}{\end{itemize}}
\newcommand{\bdf}{\begin{defin}}
\newcommand{\edf}{\end{defin}}

\newcommand{\nn}{\nonumber}
\renewcommand{\=}{&=&}
\renewcommand{\>}{&>&}
\renewcommand{\le}{\leq}
\newcommand{\+}{&+&}
\newcommand{\fr}{\frac}

\renewcommand{\r}{{\mathbb R}}
\newcommand{\br}{\bar{\mathbb R}}
\newcommand{\Z}{{\mathbb Z}}
\newcommand{\z}{{\mathbb Z}}
\newcommand{\zd}{\z^d}
\newcommand{\zz}{{\mathbb Z}}
\newcommand{\R}{{\mathbb R}}
\newcommand{\tw}{\tilde{\cal W}}
\newcommand{\E}{{\mathbb E}}
\newcommand{\C}{{\mathbb C}}
\renewcommand{\P}{{\mathbb P}}
\newcommand{\N}{{\mathbb N}}
\newcommand{\var}{{\mathbb V}}
\renewcommand{\S}{{\cal S}}
\newcommand{\T}{{\cal T}}
\newcommand{\W}{{\cal W}}
\newcommand{\X}{{\cal X}}
\newcommand{\Y}{{\cal Y}}
\newcommand{\cm}{{\cal M}}
\newcommand{\cp}{{\cal P}}
\newcommand{\h}{{\cal H}}
\newcommand{\f}{{\cal F}}
\newcommand{\cd}{{\cal D}}
\newcommand{\xt}{X_t}
\renewcommand{\ge}{g^{(\epsilon)}}
\newcommand{\xe}{y^{(\epsilon)}}
\newcommand{\ye}{y^{(\epsilon)}}
\newcommand{\bx}{{\bar y}}
\newcommand{\by}{{\bar y}}

\newcommand{\bxe}{{\bar y}^{(\epsilon)}}
\newcommand{\bye}{{\bar y}^{(\epsilon)}}
\newcommand{\bwe}{{\bar w}^{(\epsilon)}}
\newcommand{\bxz}{{\bar y}}
\newcommand{\bwz}{{\bar w}}
\newcommand{\we}{w^{(\epsilon)}}
\newcommand{\Xe}{Y^{(\epsilon)}}
\newcommand{\Ze}{Z^{(\epsilon)}}
\newcommand{\Ye}{Y^{(\epsilon)}}
\newcommand{\ydo}{Y^{(\d)}_{y_0(\d),s_0(\d)}}
\newcommand{\yo}{Y_{y_0,s_0}}
\newcommand{\tye}{{\tilde Y}^{(\epsilon)}}
\newcommand{\hy}{{\hat Y}}
\newcommand{\ve}{V^{(\epsilon)}}
\newcommand{\Ne}{N^{(\epsilon)}}
\newcommand{\ce}{c^{(\epsilon)}}
\newcommand{\cle}{c^{(\l\epsilon)}}
\newcommand{\xet}{Y^{(\epsilon)}_t}
\newcommand{\hxt}{\hat X_t}
\newcommand{\btn}{\bar\tau_n}
\newcommand{\ct}{{\cal T}}
\newcommand{\rn}{{\cal R}_n}
\newcommand{\nt}{{N}_t}
\newcommand{\lnk}{{\cal L}_{n,k}}
\newcommand{\cl}{{\cal L}}
\newcommand{\bw}{\bar{\cal W}}
\newcommand{\tc}{\tilde{\cal C}_b}
\newcommand{\hxtt}{\hat X_{\ct}}
\newcommand{\txnt}{\tilde X_{\nt}}
\newcommand{\xs}{X_s}
\newcommand{\xn}{\tilde X_n}
\newcommand{\tx}{\tilde X}
\newcommand{\hx}{\hat X}
\newcommand{\txi}{\tilde X_i}
\newcommand{\txij}{\tilde X_{i_j}}
\newcommand{\taxi}{\tau_{\txi}}
\newcommand{\txn}{\tilde X_N}
\newcommand{\xk}{X_K}
\newcommand{\ts}{\tilde S}
\newcommand{\tl}{\tilde\l}
\newcommand{\tg}{\tilde g}
\newcommand{\im}{I^-}
\newcommand{\ip}{I^+}
\newcommand{\hal}{H_\a}
\newcommand{\ba}{B_\a}
\renewcommand{\a}{\alpha}
\renewcommand{\b}{\beta}
\newcommand{\g}{\gamma}
\newcommand{\G}{\Gamma}
\renewcommand{\L}{\Lambda}
\renewcommand{\d}{\delta}
\newcommand{\D}{\Delta}
\newcommand{\e}{\epsilon}
\newcommand{\fes}{\phi^{(\epsilon)}_s}
\newcommand{\fet}{\phi^{(\epsilon)}_t}
\newcommand{\fe}{\phi^{(\epsilon)}}
\newcommand{\pset}{\psi^{(\epsilon)}_t}
\newcommand{\pse}{\psi^{(\epsilon)}}
\renewcommand{\l}{\lambda}
\newcommand{\me}{\mu^{(\epsilon)}}
\newcommand{\re}{\rho^{(\epsilon)}}
\newcommand{\tre}{\tilde{\rho}^{(\epsilon)}}
\newcommand{\nue}{\nu^{(\epsilon)}}
\newcommand{\mbe}{{\bar\mu}^{(\epsilon)}}
\newcommand{\rbe}{{\bar\rho}^{(\epsilon)}}
\newcommand{\mb}{{\bar\mu}}
\newcommand{\rb}{{\bar\rho}}
\newcommand{\mbz}{{\bar\mu}}
\newcommand{\s}{\sigma}
\renewcommand{\o}{\Pi}
\newcommand{\om}{\omega}
\newcommand{\tio}{\tilde\o}
\renewcommand{\sl}{\sigma'}
\newcommand{\si}{\s(i)}
\newcommand{\sit}{\s_t(i)}
\newcommand{\ei}{\eta(i)}
\newcommand{\eit}{\eta_t(i)}
\newcommand{\eot}{\eta_t(0)}
\newcommand{\sil}{\s'_i}
\newcommand{\sj}{\s(j)}
\newcommand{\st}{\s_t}
\newcommand{\so}{\s_0}
\newcommand{\xii}{\xi_i}
\newcommand{\xij}{\xi_j}
\newcommand{\xio}{\xi_0}
\newcommand{\ti}{\tau_i}
\newcommand{\te}{\tau^{(\epsilon)}}
\newcommand{\bt}{\bar\tau}
\newcommand{\tti}{\tilde\tau_i}
\newcommand{\tto}{\tilde\tau_0}
\newcommand{\tei}{T_i}
\newcommand{\ttei}{\tilde T_i}
\newcommand{\tes}{T_S}
\newcommand{\tao}{\tau_0}

\renewcommand{\t}{\tilde t}

\newcommand{\da}{\downarrow}
\newcommand{\ua}{\uparrow}
\newcommand{\ar}{\rightarrow}
\newcommand{\lar}{\leftrightarrow}
\newcommand{\va}{\stackrel{v}{\rightarrow}}
\newcommand{\ppa}{\stackrel{pp}{\rightarrow}}
\newcommand{\dw}{\stackrel{w}{\Rightarrow}}
\newcommand{\Va}{\stackrel{v}{\Rightarrow}}
\newcommand{\Ppa}{\stackrel{pp}{\Rightarrow}}
\newcommand{\la}{\langle}
\newcommand{\ra}{\rangle}
\newcommand{\ep}{\vspace{.5cm}}
\newcommand\sqr{\vcenter{
              \hrule height.1mm
              \hbox{\vrule width.1mm height2.2mm\kern2.18mm\vrule
width.1mm}
              \hrule height.1mm}}                  % This is a slimmer sqr.

\newcommand{\stack}[2]{\genfrac{}{}{0pt}{3}{#1}{#2}}

%%%%%%%%%%%%%%%%%%%%%%%%%%%%%%%%%%%%%%%%%%%%%%%
%%%%%%%%%%%%%%% ABSTRACT %%%%%%%%%%%%%%%%%%%%%%%%%%
%%%%%%%%%%%%%%%%%%%%%%%%%%%%%%%%%%%%%%%%%%%%%%%

\begin{abstract}

The Brownian Web (BW) is
the random network formally
consisting of the paths
of coalescing one-dimensional Brownian motions
starting from every space-time point in $\r\times\r$.
We extend the earlier work of Arratia
and of T\'oth and Werner by providing characterization and
convergence results for the BW distribution, including convergence of the
system of all coalescing random walks to the BW under diffusive space-time
scaling. We also
provide characterization and convergence results for the Double Brownian Web,
which combines the BW with its dual process of coalescing Brownian
motions moving backwards in time,  with forward and backward paths
``reflecting'' off each other.
For the BW, deterministic space-time points
are almost surely of ``type'' $(0,1)$ --- {\em zero} paths into the point from the
past and exactly
{\em one} path out of the point to the future;
we determine the Hausdorff dimension for all types that actually
occur: dimension $2$ for
type $(0,1)$,
$3/2$ for $(1,1)$ and $(0,2)$, $1$ for $(1,2)$, and $0$ for $(2,1)$ and $(0,3)$.

\end{abstract}



%%%%%%%%%%%%%%%%%%%%%%%%%%%%%%%%%%%%%%%%%%%%%%%
%%%%%%%%%%%%%% INTRODUCTION%%%%%%%%%%%%%%%%%%%%%%%%
%%%%%%%%%%%%%%%%%%%%%%%%%%%%%%%%%%%%%%%%%%%%%%%
\section{Introduction}

\setcounter{equation}{0}
\label{sec:intro}

In this paper, we present a number of results
concerning the characterization of and convergence to
a striking stochastic object called the {\em Brownian Web}
(as well as similar results for the closely related
{\em Double Brownian Web}). Several of the main results were
previously announced, with sketches of
the proofs, in~\cite{kn:FINR}.

Roughly speaking, the Brownian Web is
the collection of graphs of coalescing one-dimensional Brownian motions
(with unit diffusion constant and zero drift) starting from
all possible starting points in one plus one dimensional (continuous)
space-time. This object was originally studied more than twenty
years ago by Arratia~\cite{kn:A2}, motivated by asymptotics
of one-dimensional voter models, and then about five years ago
by T\'{o}th and Werner~\cite{kn:TW}, motivated by the problem of
constructing continuum ``self-repelling motions.'' Our
own interest in this object arose because of its relevance to
``aging'' in statistical
physics models of one-dimensional
coarsening~\cite{kn:FINS1,kn:FINS2}---which returns us to
Arratia's original context of voter models, or equivalently
coalescing random walks in one dimension. This
motivation leads to our primary concern
with weak convergence results, which in turn requires a
careful choice of space for the Brownian Web so as to
obtain useful characterization criteria for its distribution.
We continue the introduction by discussing coalescing random
walks and their scaling limits.

Let us begin by constructing random paths
in the plane, as follows.
Consider the two-dimensional lattice of all points
$(i,j)$ with $i,j$ integers and $i+j$ even.
Let a walker at spatial location $i$ at time $j$
move right or left at unit speed between times $j$ and $j+1$
if the outcome of a fair coin toss is heads ($ \Delta_{i,j} = +1 $)
or tails ($ \Delta_{i,j} = -1 $), with the coin tosses independent
for different space-time points $(i,j)$. Figure 1 depicts a
simulation of the resulting paths.

%%%%%%%%%%%%%%%%%%%%%%%%%%%%%%%%%%%%%%%%%%%%
%%%%%%%%%%%%%%  FIGURE 1 %%%%%%%%%%%%%%%%%%%%%%%%
%%%%%%%%%%%%%%%%%%%%%%%%%%%%%%%%%%%%%%%%%%%%
\begin{figure}[!ht]
\begin{center}
\includegraphics[width=10cm]{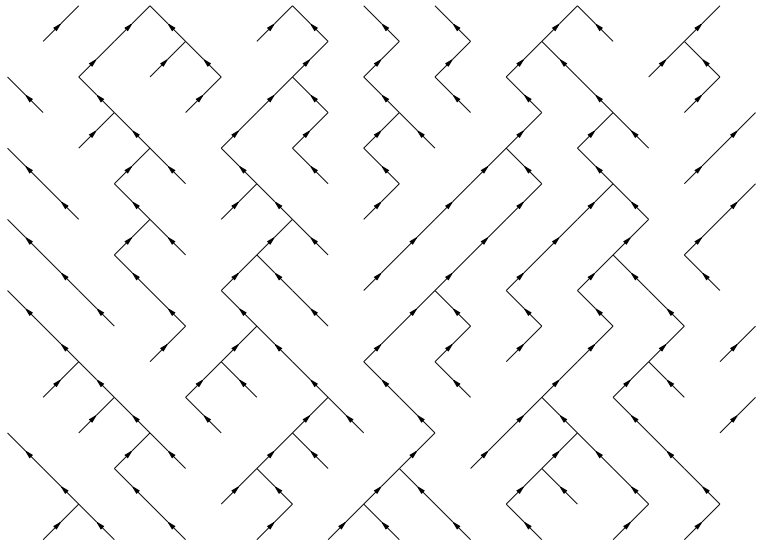}
\caption{Coalescing random walks in discrete time;
the horizontal coordinate is space and the
vertical one is time.}\label{figpath}
\label{path} 
\end{center}
\end{figure}

The path of a walker starting from $y_0$ at time $s_0$ is the graph
of a simple
symmetric one-dimensional random walk, $\yo (t)$. At integer times, $\yo (t)$
is the solution of the simple stochastic difference equation,
\begin{equation}
\label{diffeq}
Y(j+1) - Y(j) =  \Delta_{Y(j),j} , \quad Y(s_0) = y_0.
\end{equation}
Furthermore the paths of distinct walkers starting from different
$(y_0,s_0)$'s are automatically {\em coalescing} --- i.e., they are
independent of each other
until they coalesce (i.e., become identical) upon meeting
at some space-time point.

If the increments $\Delta_{i,j}$ remain i.i.d., but
take values besides $\pm1$ (e.g., $\pm3$), then one obtains
non-simple random walks whose paths can cross each other in
space-time, although they still coalesce when they land
on the same space-time lattice site. Such systems with
crossing paths will be discussed in Section~\ref{sec:conv}.

After rescaling to spatial steps of size $\d$ and time steps
of size $\d^2$, a single rescaled random walk (say, starting from
$0$ at time $0$) $Y_{0,0}^{(\d)}(t) = \d Y_{0,0}(\d^{-2}t)$
converges as $\d \to 0$ to a standard Brownian motion $B(t)$.
 That is,
by the Donsker invariance principle~\cite{kn:D},
the distribution of $Y_{0,0}^{(\d)}$ on the space of continuous paths
converges as $\d\to0$ to standard Wiener measure.

The invariance principle is also valid for continuous time random walks,
where the move from $i$ to $i\pm1$ takes an exponentially distributed
time .
In continuous time, coalescing random walks are at the
heart of Harris's graphical representation
of the (one-dimensional) voter model~\cite{kn:H} and their scaling limits
arise naturally in the physical context of (one-dimensional) aging
(see, e.g.,~\cite{kn:FINS1,kn:FINS2}).
Of course, finitely many rescaled coalescing walks
in discrete or continuous time (with rescaled space-time starting points)
converge in distribution to  finitely many coalescing
Brownian motions.
In this paper, we present results
concerning the convergence
in distribution of the complete
collection of the rescaled coalescing walks from
{\em all} the starting points.

Our results are in two main parts:
\begin{itemize}
\item[(1)] characterization (and construction) of the limiting object,
the standard {\em Brownian Web (BW)}, and
\item[(2)] convergence criteria, which are applied, in this paper, to
coalescing random walks.
\end{itemize}

 As a cautionary remark, we point out that the scaling limit
motivating our convergence results does not belong to the
realm of hydrodynamic limits of particle systems but
rather to the realm of invariance principles.

A key ingredient of the characterization and construction
(see, e.g., Theorem~\ref{teo:char}) 
is the choice of a space for the Brownian web; this is the
BW analogue of the space of continuous paths for Brownian motion.
The convergence criteria and application (see, e.g.,
Theorems~\ref{n-c-conv} and~\ref{teo:convrw} below)
are the BW analogues of Donsker's invariance principle.
Like Brownian motion itself,
we expect that the Brownian web
and its variants (see, e.g., Sec.~\ref{sec:double}) will be quite
ubiquitous as scaling limits, well beyond the context of
coalescing random walks
(and our sufficient conditions for convergence);
one situation where this occurs is for two-dimensional
``Poisson webs''~\cite{kn:FFW}. Another example is in the area
of river basin modelling. In~\cite{kn:S}, coalescing random walks
were proposed as a model of a drainage network. Some of the questions
about scaling in such models may find answers in the context of their
scaling limits. For more on the random walk and other models for river
basins, see~\cite{kn:RR}.

Much of the construction of the Brownian web
(but without characterization or convergence results)
was already done in the groundbreaking work of
Arratia~\cite{kn:A2,kn:A0} (see also~\cite{kn:A1,kn:H2})
and then in  the work of T\'oth and
Werner~\cite{kn:TW} (see also~\cite{kn:STW}
and~\cite{kn:T};  in the latter reference, the Brownian web is introduced
in relation to {\em Black Noise}).
Arratia, T\'oth and Werner all recognized that in the
limit $\d\to0$ there would be (nondeterministic) space-time points
$(x,t)$ starting from which
there are multiple limit paths and they provided
various conventions (e.g., semicontinuity in $x$) to avoid such
multiplicity. Our main contribution vis-a-vis
construction is to accept the intrinsic
nonuniqueness by choosing an appropriate metric space in which the BW takes
its values. Roughly speaking, instead of using some convention to obtain
a process that is a {\em single-valued} mapping from
each space-time starting
point to a single path from that starting point, we allow
{\em multi-valued} mappings; more accurately, our BW value is the collection
of {\em all} paths from all starting points. This choice of space is
very much in the spirit of earlier
work~\cite{kn:A,kn:AB,kn:ABNW} on spatial scaling limits of
critical percolation models and
spanning trees, but modified for our particular space-time setting;
the directed (in time) nature of our paths considerably simplifies
the topological setting compared to~\cite{kn:A,kn:AB,kn:ABNW}.

The Donsker invariance principle implies that the distribution
of any continuous (in the sup-norm metric) functional of
$Y_{0,0}^{(\d)}$ converges to that for Brownian motion. The
classic example of such a functional is the random walk
maximum, $\sup_{0 \le t \le 1}Y_{0,0}^{(\d)}(t)$.
An analogous example
for coalescing random walks
is the maximum over all rescaled walks
starting at (or passing through) some vertical (time-like) interval,
i.e., the maximum value (for times $t \in [s,1]$) over walks
touching any space-time point of the form
$(0, s)$ for some $s \in [0,1]$.
In this case, the functional is not quite continuous
for our choice of metric space, but it is continuous almost everywhere
(with respect to the Brownian web measure), which is sufficient.

The rest of the paper is organized as follows. In Section~\ref{sec:results},
we state two of the main results of the paper:
Theorem~\ref{teo:char} is a characterization of the
Brownian Web and Theorem~\ref{n-c-conv}
is a  convergence theorem for the
important special case where, even before taking a limit, all paths are
noncrossing. Section~\ref{sec:char} contains the proof
of the initial characterization result, Theorem~\ref{teo:char},
as well as an alternative characterization, Theorem~\ref{teo:charm}
in which a kind of separability condition is replaced
by a minimality condition. In Section~\ref{sec:count}, we present
other characterization
results (Theorems~\ref{teo:char1} and~\ref{teo:char2}) based on certain
counting  random variables,
which will be needed for the derivation of our
main convergence results. Section~\ref{sec:double} concerns the
construction and characterization of the {\em Double Brownian Web\/}
which is the union of the Brownian Web with its dual system of
coalescing Brownian motions going backwards in time. Our analysis
of the Double Brownian Web relies heavily on the work of~\cite{kn:STW}.
Section~\ref{sec:double} also contains an analysis of the typology of
generic and exceptional space-time points for the Brownian web, as
in~\cite{kn:TW}, and of the Hausdorff dimension of the various types.
In Section~\ref{sec:conv}, we extend our convergence results to
cover the case of crossing paths; the proof of the noncrossing
result, Theorem~\ref{n-c-conv}, is given here as a corollary
of the more general result. In  Section~\ref{sec:rw-conv}
we apply our (non-crossing) convergence results to the case of
coalescing random walks. There are three appendices,
the first covers issues of mesurability, the second issues of compactness
and tightness, and the third gives a Hausdorff dimension result about
Brownian motion graphs that is used in Section~\ref{sec:double}.

%%%%%%%%%%%%%%%%%%%%%%%%%%%%%%%%%%%%%%%%%%%%
%%%%%%%%%%%%%%%% SECTION 2 %%%%%%%%%%%%%%%%%%%%%%
%%%%%%%%%%%%%%%%%%%%%%%%%%%%%%%%%%%%%%%%%%%%
%%%%%%%%%%%%%%%%%%%%%%%%%%%%%%%%%%%%%%%%%%%%

\section{Some Main Results}
\setcounter{equation}{0}
\label{sec:results}

We begin by defining three metric spaces: $(\br^2,\rho)$, $(\o,d)$ and
$(\h,d_\h)$.
The elements of the three spaces are respectively: points
in space-time, paths with specified starting points in space-time and
collections of paths with specified starting points.
The BW will be an $(\h,\f_\h)$-valued random variable,
where $\f_\h$ is the Borel $\s$-field  associated to the metric $d_\h$.

$(\br^2,\rho)$ is the completion (or compactification) of $\R^2$ under
a metric $\rho$ given below (see~(\ref{rho}) and~(\ref{compactify})).
$\br^2$ may be thought as the set of $(x,t)$ in
$[-\infty,\infty]\times[-\infty,\infty]$ with all points of the form
$(x,-\infty)$ identified (and similarly for $(x,\infty)$).

For $t_0\in[-\infty,\infty]$, let $C[t_0]$ denote the set of functions
$f$ from $[t_0,\infty]$ to $[-\infty,\infty]$ such that
$(1+|t|)^{-1} \tanh (f(t))$ is continuous. Then define
\begin{equation}
%\label{omega}
\o=\bigcup_{t_0\in[-\infty,\infty]}C[t_0]\times\{t_0\},
\end{equation}
where $(f,t_0)\in\o$ represents a path in $\br^2$ starting at
$(f(t_0),t_0)$.
For$(f,t_0)$ in $\o$, we denote by $\hat f$ the function that extends $f$
to all
$[-\infty,\infty]$ by setting it equal to $f(t_0)$ for $t<t_0$. Then we
define a suitable distance $d$ (see~(\ref{d})) such that
$(\Pi,d)$ is a complete separable metric space.

Finally, $\h$ denotes the set of compact
subsets of $(\Pi,d)$, with $d_\h$ the
induced Hausdorff metric (see~(\ref{dh}))
and $\f_\h$ the Borel $\sigma$-algebra.
$(\h,d_\h)$ is also a complete separable metric space. For
an $(\h,\f_\h)$-valued random variable $\bar \W$
(or its
distribution $\mu$), we define the {\em finite-dimensional
distributions}
of $\bar \W$ as the induced probability measures
$\mu_{(x_1,t_1;\ldots;x_n,t_n)}$
on the subsets of paths starting from any
finite deterministic set of points
$(x_1,t_1),\ldots,(x_n,t_n)$ in $\R^2$.
There are several ways in which the Brownian web can be
characterized; they differ from each other primarily
in the type of extra condition required
beyond the finite-dimensional distributions (which are those
of coalescing Brownian motions).
In the next theorem,  the extra condition is a type of Doob separability
property (see, e.g., Chap.~3 of \cite{kn:V}).
Variants are stated later either using a minimality property
(Theorem~\ref{teo:charm}) or a
counting  random variable (Theorems~\ref{teo:char1} and
\ref{teo:char2}). Theorem~\ref{teo:char2} is the one most directly
suited to the convergence results of Section~\ref{sec:conv}.

The events and random variables appearing in the next two theorems
are  $(\h,\f_\h)$-measurable. 
  This claim follows straightforwardly from Proposition A.2 in
 Appendix \ref{app:mea}.
The proofs of the two theorems are given, respectively, in
Sections~\ref{sec:char} and~\ref{sec:conv}.

\bteo
\label{teo:char}
There is an \( ({\cal H},{\cal F}_{{\cal H}}) \)-valued random variable \(
\bar{\W} \)
whose distribution is uniquely determined by the following three
properties.
\begin{itemize}
         \item[(o)]  from any deterministic point \( (x,t) \) in
$\r^{2}$,
         there is almost surely a unique path \( {W}_{x,t} \)
starting from \( (x,t) \).

         \item[(i)]  for any deterministic \( n, (x_{1}, t_{1}), \ldots,
         (x_{n}, t_{n}) \), the joint distribution of \(
         {W}_{x_{1},t_{1}}, \ldots, {W}_{x_{n},t_{n}} \) is that
         of coalescing Brownian motions (with unit diffusion constant),
and

         \item[(ii)]  for any deterministic, dense countable subset  \(
{\cal
         D} \) of \( \r^{2} \), almost surely, \( \bar{\W} \) is the
closure in
         \( ({\cal H}, d_{{\cal H}}) \) of \( \{ {W}_{x,t}: (x,t)\in
         {\cal D} \}. \)
\end{itemize}
\eteo

\brm
One can choose a
{\em single} dense countable ${\cal D}_{0}$ and in $(o)$, $(i)$ and $(ii)$
restrict to space-time starting points from that ${\cal D}_{0}$.
Different characterization theorems for the Brownian web with alternatives for $(ii)$
are given in Sections~\ref{sec:char} and~\ref{sec:count}.
We note that there are
natural \( ({\cal H},{\cal F}_{{\cal H}}) \)-valued random variables satisfying
$(o)$ and $(i)$ but not $(ii)$.  An instance of such a random variable
(closely related to the double
Brownian web) will be studied elsewhere,
and shown to arise as the scaling limit
of stochastic flows, extending earlier work of Piterbarg~\cite{kn:P}.
\erm

The next theorem is restricted to noncrossing
processes, but a somewhat more general result is given in Section~\ref{sec:conv}.

\bdf
\label{eta}
For $t>0,\,t_0,a,b\in\r,\,a<b$, let $\eta(t_0,t;a,b)$ be the number
of {\em distinct} points in
$\R\times\{t_0+t\}$ that are touched by paths in $\bar \W$ which also
touch
some point in $[a,b]\times\{t_0\}$. Let also
$\hat{\eta}(t_{0}, t; a,b)=\eta(t_{0}, t; a,b)-1.$
\edf

We note that by duality arguments (see Remark~\ref{rmk:dual}), it can be shown
that for deterministic $t_0,t,a,b$, this $\hat\eta$ is equidistributed with the
number of distinct points in $[a,b]\times\{t_0+t\}$ that are touched by paths
in $\bar\W$ which also touch $\r\times\{t_0\}$.

\bteo
\label{n-c-conv}
Suppose \( \X_{1}, \X_{2}, \dots \) are
\(({\cal H}, {\cal F}_{{\cal H}}) \)-valued random variables with
noncrossing
paths. If, in addition, the following three conditions are valid, then
the distribution
$\mu_{n}$ of $\X_{n}$ converges to the distribution $\mu_{\bar{W}}$ of the
standard Brownian web.

\begin{itemize}
         \item[$(I_{1})$]  There exist such $\theta^y_n\in \X_n$
                          satisfying:
                         for any deterministic $y_1,\ldots,y_m\in\cd$,
                         ${\theta^{y_1}_n,\ldots,\theta^{y_m}_n}$
                        converge
                      in distribution as $n\to \infty$ to coalescing Brownian
                      motions (with unit diffusion constant) starting at
                      $y_1,\ldots,y_m$.
         \item[$(B_{1})$]
               $\forall t>0,\,\,{\displaystyle
               \limsup_{n \to \infty}\sup_{(a,t_0)\in\R^2}}
                         \mu_n(\hat{\eta}(t_0,t;a,a+\e)\geq1)\to0
                          \hbox{ as } \e\to0+$;
         \item[$(B_{2})$]
          $\forall t>0,\,\,\e^{-1}{\displaystyle
          \limsup_{n \to \infty}\sup_{(a,t_0)\in\R^2}}
                         \mu_n(\hat{\eta}(t_0,t;a,a+\e)\geq2)\to0
                      \hbox{ as } \e\to0+$.
\end{itemize}
\eteo

Convergence of coalescing random walks (in discrete and
continuous time) (see Theorem~\ref{teo:convrw})
is obtained
as a
corollary to Theorem~\ref{n-c-conv}.

%%%%%%%%%%%%%%%%%%%%%%%%%%%%%%%%%%%%%%%%%%%%%%
%%%%%%%%%%%%%%%% SECTION 3 %%%%%%%%%%%%%%%%%%%%%%%%
%%%%%%%%%%%%%%%%%%%%%%%%%%%%%%%%%%%%%%%%%%%%%%
%%%%%%%%%%%%%%%%%%%%%%%%%%%%%%%%%%%%%%%%%%%%%%

\section{Construction and initial characterizations}

\setcounter{equation}{0}
\label{sec:char}

In this section we give a complete construction of the Brownian web
which yields a proof of Theorem \ref{teo:char}. Then we give in
Theorem~\ref{teo:charm} a somewhat
different characterization of the BW distribution.

Let $(\Omega,{\cal F}, \P)$ be a probability space where an
i.i.d.~family of standard Brownian motions $(B_j)_{j\geq1}$
is defined. Let ${\cal D}=\{(x_j,t_j),\,j\geq1\}$ be a countable dense
set in $\R^2$.
Let $W_j$ be a Brownian
path starting at position $x_j$ at time $t_j$. More precisely,
\begin{equation}
\label{ibp}
W_j(t)=x_j+B_j(t-t_j),\, \, t\geq t_j.
\end{equation}

We now construct coalescing Brownian paths out of the family of paths
$(W_j)_{j\geq1}$ by specifying coalescing rules. When two paths meet
for the first time, they coalesce into a single path, which is that of
the Brownian motion with the lower label. We denote the coalescing
Brownian
paths by $\tilde W_j,\,j\geq1$.
Notice the strong Markov property of Brownian motion allows for a lot
of freedom in giving a coalescing rule. Any rule, even non local, that
does not depend on the realization of the $(W_j)$'s {\em after} the
time of coalescence will yeld the same object in distribution

We now give a more careful inductive definition of the coalescing
paths.
First we set
\begin{equation}
\label{cbp1}
\tilde W_1=W_1.
\end{equation}
For $j\geq2$, let $\tilde W_j$ be the mapping from $[t_j,\infty)$
to $\R$ defined as follows. Let
\beqn
\label{cbp2}
\tau_j&=&\inf\{t\geq t_j:\,W_j(t)=\tilde W_i(t)
\mbox{ for some } 1\leq i\leq j-1\},\\
\label{cbp3}
I_j&=&\min\{1\leq i\leq j-1:\,W_j(\tau_j)=\tilde W_i(\tau_j)\},
\eeqn
and define
\beqn
\nonumber
\tilde W_j(t)&=&W_j(t),\mbox{ if }t_j\leq t\leq\tau_j\\
\label{cbp4}
                &=& \tilde W_{I_j}(t),\mbox{ if }t>\tau_j.
\eeqn
We then define the
Brownian web skeleton ${\cal W}({\cal D})$
with {\em starting set} ${\cal D}$ by
\begin{eqnarray}
{\cal W}_{k} &=&{\cal W}_{k}({\cal D})=\{\tilde W_j: \, 1\leq j\leq k \}
\\
{\cal W} &=& {\cal W}({\cal D})=\bigcup_{k} {\cal W}_{k}
\end{eqnarray}

$(\br^2,\rho)$ is the completion (or compactification) of $\R^2$ under
the
metric $\rho$, where
\begin{equation}
\label{rho}
\rho((x_1,t_1),(x_2,t_2))=
\left|\frac{\tanh(x_1)}{1+|t_1|}-\frac{\tanh(x_2)}{1+|t_2|}\right|
\vee|\tanh(t_1)-\tanh(t_2)|.
\end{equation}
$\br^2$ may be thought as
the image of $[-\infty,\infty]\times[-\infty,\infty]$
under the mapping
\begin{equation}
\label{compactify}
(x,t)\leadsto(\Phi(x,t),\Psi(t))
\equiv\left(\frac{\tanh(x)}{1+|t|},\tanh(t)\right).
\end{equation}

For $t_0\in[-\infty,\infty]$, let $C[t_0]$ denote the set of functions
$f$ from $[t_0,\infty]$ to $[-\infty,\infty]$ such that $\Phi(f(t),t)$
is continuous. Then define
\begin{equation}
%\label{omega}
\o=\bigcup_{t_0\in[-\infty,\infty]}C[t_0]\times\{t_0\},
\end{equation}
where $(f,t_0)\in\o$ represents a path in $\br^2$ starting at
$(f(t_0),t_0)$.
For$(f,t_0)$ in $\o$, we denote by $\hat f$ the function that extends $f$
to all
$[-\infty,\infty]$ by setting it equal to $f(t_0)$ for $t<t_0$. Then we
take
\begin{equation}
\label{d}
d((f_1,t_1),(f_2,t_2))=
(\sup_t|\Phi(\hat{f_1}(t),t)-\Phi(\hat{f_2}(t),t)|)
\vee|\Psi(t_1)-\Psi(t_2)|.
\end{equation}
$(\o,d)$ is a complete separable metric space.

Let now $\h$ denote the set of compact
subsets of $(\o,d)$, with $d_\h$ the induced Hausdorff metric, i.e.,
\begin{equation}
\label{dh}
d_\h(K_1,K_2)=\sup_{g_1\in K_1}\inf_{g_2\in K_2}d(g_1,g_2)\vee
                 \sup_{g_2\in K_2}\inf_{g_1\in K_1}d(g_1,g_2).
\end{equation}
$(\h,d_\h)$ is also a complete separable metric space.

\bdf
$\bar{{\cal W}} ({\cal D})$ is the closure in $(\Pi,d)$ of ${\cal
W}({\cal D})$.
\edf

\bprop
\label{prop:findim}
\( \bar{{\cal W}} ({\cal D}) \) satisfies properties (o) and (i)
of Theorem \ref{teo:char}; i.e., its finite dimensional distributions
(whether from points in \( {\cal D} \) or not) are those of coalescing
Brownian motions.
\eprop
{\bf Proof }
We will prove the proposition in the special case of the
distribution of path(s) from the single point at the origin $(0,0)$.
The proof easily  extends to the case of finitely many space-time points.
Let  $\{(x_n^r,t_n^r)\}_{n=1}^\infty$ and $\{(x_n^l,t_n^l)\}_{n=1}^\infty$ be
two
sequences of points from ${\cal D}$ satisfying the following conditions.
  There exist positive constants $c_1,c_2 \leq 1/2$ such that
\beqnn
&&-c_1/{n}^2 < x_n^l <0 <  x_n^r < c_2/{n}^2,\\
&&-|x_n^l|^3 < t_n^l < 0, -|x_n^r|^3 < t_n^r < 0.
\eeqnn

\noindent Let $\tilde W_{n,l}$ and $ \tilde W_{n,r}$ be the
coalescing Brownian motions
starting from $(x_n^l,t_n^l)$ and $(x_n^r,t_n^r)$ respectively and let
$B(t)$ be a standard
Brownian motion (which for convenience we take as defined on
$(\Omega,{\cal F}, \P)$,
e.g. by letting $B = B_1$).

Using the reflection principle we obtain
\beqn
\nonumber
&&\P\left(\max_{t_n^l \leq s \leq 0} \tilde W_{n,l}(s) \geq 0\right) =  2 \P \left(
B({|t_n^l|}) \geq |x_n^l|\right) \\
\nonumber
& = & 2 \P \left( B(1) \geq |x_n^l|/|t_n^l|^{1/2}\right) 
\leq  
2 \P \left(B(1) \geq 1 / {|x_n^l|}^{1/2}\right) \\
\label{lcr}
& \leq  & 2 \P \left(B(1) \geq \sqrt 2 \,n\right)  \leq  e^{-n^2}.
\eeqn

\noindent Proceeding along the same lines we can prove that
\begin{equation}
\label{rcr}
\P\left(\min_{t_n^r \leq s \leq  0} \tilde W_{n,r}(s) \leq 0\right) \leq
e^{-n^2} .
\end{equation}

\noindent Using scaling properties of Brownian motion,  we obtain
\beqn
\nonumber
&&\P\left(\tilde W_{n,l}(0) \leq  -2 c_1/n^2\right)  \leq  \P\left(B({|t_n^l|})\geq c_1/n^2\right) \\
\nonumber
& = & \P\left(B(1) \geq \frac{c_1}{n^2} \frac{1}{|t_n^l|^{1/2}}\right) 
\leq   \P\left(B(1) \geq \frac{c_1}{n^2} \frac{n^3}{c_1^{3/2}}\right)\\
\label{zero}
& = & \P\left(B(1) \geq n/c_1^{1/2}\right)
\leq  \frac12\,  e^{-n^2/(2c_1)} \leq e^{- n^2}.
\eeqn
\noindent Similarly,
\[
\P\left(\tilde W_{n,r}(0) \geq 2 c_2/n^2\right) \leq  e^{-n^2}.
\]

Now, let  $\tau_n = \max \{t_n^l,t_n^r\}$.
Then, for all $n\geq1$ we have the
following inequalities and equalities
(where we use the reflection principle in
one place).
\begin{eqnarray}
\nonumber
&&\P\left(\inf_{\tau_n\leq s\leq 1/n}|\tilde W_{n,l}(s)-\tilde W_{n,r}(s)|>0\right)\\
\nonumber
&\leq &
\P\left(\tilde W_{n,l}(0) \leq  -2 c_1/n^2\right) + \P\left(\tilde W_{n,r}(0) \geq 2 c_2/n^2\right)\\
\nonumber
&+&\P\left(\inf_{0\leq s\leq 1/n}|\tilde W_{n,l}(s)-\tilde W_{n,r}(s)|>0,
\tilde W_{n,l}(0)>-\frac{2c_1}{n^2},
\tilde W_{n,r}(0)<\frac{2c_2}{n^2}\right)\\
\nonumber
&\leq&2e^{-n^2}+1-\P\left(\inf_{0\leq s\leq 1/n}
\left[\sqrt2B(s)+2(c_1+c_2)/n^2\right]>0\right)\\
\nonumber
&=&2e^{-n^2}+1-\P\left(\sup_{0\leq s'\leq 2/n}
B(s')<2(c_1+c_2)/n^2\right)\\
\nonumber
&=&2e^{-n^2}+1-2\P\left(B(2/n)>2(c_1+c_2)/n^2\right)\\
\nonumber
& = &  2  e^{-n^2} +1 - 2 \P\left( B(1) > \frac{2(c_1+c_2)/n^2}{\sqrt{2/n}}\right)  \\
\label{lrmeetT}
& \leq&C/n^{3/2},
\end{eqnarray}
for some constant $C$.
We also observe that
\beqn
\nonumber
\P\left( \max_{t_n^l \leq s \leq 1/n}\tilde W_{n,l}(s) \geq n^{-1/4}\right) &
= &
2\P\left(B(-t_n^l+1/n) \geq {n^{-1/4}} +|x_n^l|\right) \\
\nonumber
& \leq& 2 \P \left(B(1) \geq n^{-1/4}/\sqrt{-t_n^l+1/n} \right)\\
\label{lmeetX}
& \leq&  e^{-C'n^{1/2}}
\eeqn
for some $C'$. Similarly,
\begin{equation}
\label{rmeetX}
\P\left( \inf_{t_n^l < s < 1/n}\tilde W_{n,r}(s) \leq - n^{-1/4}\right)
\leq e^{-C'n^{1/2}}.
\end{equation}

Let $\tilde W_{n,l}(0) = l_n$ and $\tilde W_{n,r}(0) = r_n$.
Put $S_n = \inf\{s: \tilde W_{n,l}(s) = \tilde W_{n,r}(s)\}$ and
let
$M_n = \tilde W_{n,l}(S_n)$.
Using the Borel-Cantelli lemma, from (\ref{lcr} ),
(\ref{rcr}),  (\ref{lrmeetT}), (\ref{lmeetX}) and (\ref{rmeetX})
   we have that
almost surely for all but finitely many $n$'s, $l_n < 0 \hbox{ and } r_n
 >0$,
   $S_n <1/n$, and
$|M_n| <n^{-1/4}$. This implies the following: Let $\Delta_n$ be
the interior of the ``triangular" space-time region with vertices
$(x_n^l,t_n^l),(x_n^r,t_n^r),(M_n,S_n)$, base $[x_n^l,x_n^r]$
and sides $\tilde W_{n,l},\tilde W_{n,r}$. There exists an almost
surely finite random
variable $N$ such that for all $k > N$
any coalescing Brownian motion $\tilde W_i $ starting from
point $(x_i,t_i)
\in {\cal D}$ and $(x_i,t_i) \in \Delta_k$, $\tilde W_i(s) =\tilde
W_{k,l}(s) =\tilde W_{k,r}(s)
$ for all $s > S_k$. Since $(M_n,S_n)$ converges to $(0,0)$ almost
surely we have, for every sequence
$((x_i,t_i))_i$ in $ {\cal D}$  converging to $(0,0)$, that the
coalescing Brownian motions
starting from $(x_i,t_i)$ converge (in the path space metric) to a
unique coalescing
Brownian motion starting from $(0,0)$, independent of the sequence.
This proves Proposition \ref{prop:findim}.

The next result is contained in Proposition~\ref{prop: compact2}.
\bprop
\label{prop:bwcpt}
\( \bar{{\cal W}} ({\cal D}) \) is almost surely a compact subset of
\( (\Pi, d) \).
\eprop

\brm
\label{rmk:lim}
Almost surely,
$\bar{{\cal W}} ({\cal D})=\lim_{k\to\infty}{\cal W}_k({\cal D})$,
where the limit
is taken in $\cal H$.
\erm

\brm
\label{rmk:holder}
It can be shown by the methods discussed in Remark~\ref{holder} of the appendix that,
almost surely, all paths in \(\bar{{\cal W}}({\cal D})\) are
H\"{o}lder continuous with exponent $\alpha$, for any $\alpha<1/2$.
\erm

\bprop
\label{prop:ind}
The distribution of \( \bar{{\cal W}} ({\cal D}) \)
does not depend on \( {\cal D} \) (including its order).
Furthermore, \( \bar{{\cal W}} ({\cal D}) \) satisfies property $(ii)$
of Theorem~\ref{teo:char}.
\eprop

\noindent{\bf Proof }  Given two countable dense sets in $\r^2$,
${\cal D}_1$ and ${\cal D}_2$, we will couple versions of
$\bar{{\cal W}} ({\cal D}_1)$ and $\bar{{\cal W}} ({\cal D}_2)$
and show that they are equal.
Let ${\cal W}_k({\cal D}_2)$ be as in the
construction
of the Brownian web with starting set ${\cal D}_2$ and let ${\cal
W}'_k({\cal D}_2)$
be the paths of $\bar{\cal W}({\cal D}_1)$ starting from the first $k$
elements
of ${\cal D}_2$. By Proposition~\ref{prop:findim}, ${\cal W}_k({\cal
D}_2)$ and
${\cal W}'_k({\cal D}_2)$ have the same distribution (namely, that of
coalescing
Brownian motions starting from the first $k$ elements of ${\cal D}_2$).
 From Proposition~\ref{prop: D_1-D_2},
$$\bar{\cal W}({\cal D}_2)\sim \bar{{\cal W}'}({\cal
D}_2) \equiv \overline{\cup_k {\cal W}'_k({\cal D}_2)}
\subseteq \bar{\cal W}({\cal D}_1).$$
By using property $(o)$ for $\bar{\cal W}({\cal D}_1)$ (more
generally, by Proposition~\ref{prop:findim}), we claim that the paths
in $\bar{\cal W}({\cal D}_1)$
starting from ${\cal D}_1$ belong to $\bar{\cal W}'({\cal D}_2)$ (
 i.e., can be approximated arbitrarily well by paths in
${\cal W}'_k({\cal D}_2)$ for sufficiently large $k$).
To see this, let $(x^1,t^1)$ be some spacetime point in ${\cal D}_1$
and let $(x^2 _j, t^2 _j)$ be a sequence in ${\cal D}_2$ converging
to $(x^1,t^1)$; then the unique (by property $(o)$) paths
$W[x^2 _j, t^2 _j]$ in $\bar{\cal W}({\cal D}_1)$ starting
from $(x^2 _j, t^2 _j)$ converge (by property $(o)$) to $W[x^1,t^1]$.
Our claim follows since $W[x^2 _j, t^2 _j]$ is in ${\cal W}'_k({\cal D}_2)$
for a sufficiently large $k$. It then follows
that $\it{all}$ paths in $\bar{\cal W}({\cal D}_1)$ can be so approximated
and hence $\bar{\cal W}({\cal D}_1) \subseteq \bar{\cal W}'({\cal D}_2)$.
Thus $\bar{\cal W}({\cal D}_1) = \bar{\cal W}'({\cal D}_2) \sim
\bar{\cal W}({\cal D}_2)$. This proves both the distributional independence
on ${\cal D}$ of $\bar{\cal W}({\cal D})$ and the validity of property $(ii)$.

\bigskip

\noindent{\bf Proof of Theorem \ref{teo:char} }
The existence part of the theorem is immediate from Propositions~\ref{prop:findim}
and~\ref{prop:ind}. The uniqueness may be argued as follows.

Take any
\(({\cal H}, {\cal F}_{{\cal H}}) \)-valued random variable $\X$ with properties
$(o), (i)$ and $(ii)$, and fix a deterministic countable dense subset
${\cal D}=\{(x_1,t_1),(x_2,t_2),\ldots\}$ of $\r^2$. Let $\bar\W=\bar{\cal W}({\cal D})$,
the version of the
Brownian web constructed at the beginning of this section of the paper, using $\cal D$
(but whose distribution
does not depend on $\cal D$, by Proposition~\ref{prop:ind}).
 From $(i)$ and Proposition~\ref{prop:findim}, $\X_n:=\{X_{x_i,t_i}:\,i=1,\ldots,n\}$, where
$X_{x,t}$ is the path of $\X$ starting at $(x,t)$ (almost surely unique by $(o)$),
is equidistributed with $\W_n:=\{W_{x_i,t_i}:\,i=1,\ldots,n\}$ for all $n\geq1$.
Now $(ii)$ and Remark~\ref{rmk:lim} imply or say that $\X=\lim_{n\to\infty}\X_n$ and
$\bar\W=\lim_{n\to\infty}\W_n$. So, they have the same distribution and the proof is complete.

The next theorem provides an alternative characterization to
Theorem \ref{teo:char}. Other characterizations that will be used for
our
convergence results, are presented in Sec. 3.

\bdf
(Stochastic ordering) \( \mu_{1} << \mu_{2} \) if for $g$ any bounded
measurable function on \(({\cal H}, {\cal F}_{{\cal H}}) \) that is
{\em increasing} (i.e., \( g(K) \leq g(K') \) when $K \subseteq K'$),
    \( \int \! g d\mu_{1} \leq \int \! g d\mu_{2} \).
\edf

\bteo
\label{teo:charm}
There is an $({\cal H}, {\cal F}_{{\cal H}})$-valued
random variable \(
\bar{\W} \)
whose distribution is uniquely determined by properties $(o),(i)$ of
Theorem~\ref{teo:char} and
\begin{itemize}
       \item[$(ii')$] if $\W^{*}$ is any other \(({\cal H}, {\cal F}_{{\cal
H}})
       \) - valued random variable satisfying (o) and (i), then
       \( \mu_{\bar{\W}} << \mu_{{\W}^{*}} \).
\end{itemize}
\eteo

\noindent{\bf Proof of Theorem \ref{teo:charm} }
The existence part of the theorem is immediate from Theorem~\ref{teo:char}.
The uniqueness may be argued as follows. Fix a deterministic countable dense subset
${\cal D}=\{(x_1,t_1),(x_2,t_2),\ldots\}$ of $\r^2$, let $\W^\ast$ be the Brownian web
constructed from $\cal D$ and let $W^\ast_{x,t}$ for $(x,t)\in{\cal D}$ be the path
from $(x,t)$ in $\W^\ast$. Since $\bar\W\supset\overline{\{W_{x,t},\,(x,t)\in{\cal D}\}}$,
where $W_{x,t}$ is the path from $(x,t)$ in $\bar\W$,
and, from $(o)$ and $(i)$, $\W_n:=\{W_{x_i,t_i},\,i=1,\ldots,n\}$ has the same distribution
as $\W_n^\ast:=\{W_{x_i,t_i}^\ast,\,i=1,\ldots,n\}$ for all $n\geq1$, we have that
\( \mu_{{\W}^{*}} << \mu_{{\bar\W}} \). This and $(ii')$ imply that
\( \mu_{{\W}^{*}} = \mu_{{\bar\W}} \), and the proof is complete.

%%%%%%%%%%%%%%%%%%%%%%%%%%%%%%%%%%%%%%%%%%%%%
%%%%%%%%%%%%%%% SECTION 4 %%%%%%%%%%%%%%%%%%%%%%%%
%%%%%%%%%%%%%%%%%%%%%%%%%%%%%%%%%%%%%%%%%%%%%
%%%%%%%%%%%%%%%%%%%%%%%%%%%%%%%%%%%%%%%%%%%%%

\section{Characterization via counting}
\setcounter{equation}{0}
\label{sec:count}

In this section, we give other characterizations of the Brownian web
that will be used for our convergence theorem.
They will be given in terms of the
counting  random variables $\eta$ and
$\hat\eta$
defined in Definition~\ref{eta}.
We begin with some properties of the Brownian web as constructed in
Section 1.

\bprop
\label{prop:eta}
For a Brownian skeleton ${\cal W}({\cal D})$, the corresponding
counting  random variable
$\hat\eta_{{\cal D}}=\hat\eta_{{\cal D}}(t_0,t;[a,b])$ satisfies
\beqn
\label{eq:fkg}
\P(\hat\eta_{{\cal D}}\geq k)
&\leq&\P(\hat\eta_{{\cal D}}\geq k-1)\P(\hat\eta_{{\cal D}}\geq 1)\\
&\leq&(\P(\hat\eta_{{\cal D}}\geq 1))^k=(\Theta(b-a,t))^k, \label{eta_D}
\eeqn
where $\Theta(b-a,t)$ is the
probability that two independent Brownian motions
starting at a distance $b-a$ apart at
time zero will not have met by time $t$ (which
itself can be expressed in terms of a single Brownian motion). Thus,
$\hat\eta_{{\cal D}}$
is almost surely finite and $\E(\hat\eta_{{\cal D}})<\infty$.
\eprop

\noindent{\bf Proof }
We note that it is sufficient to prove the inequalities for
${\cal W}_n({\cal D})= \{ \tilde{W_j} : 1\leq j \leq n\}$ for
all $n$. Morevover for a given $n$ if we prove the result for
$n$ coalescing random walks, then the result for $n$ Brownian
motions will follow from the scaling limit of the random walks.
Therefore we now consider $n$ discrete time simple symmetric coalescing
random walks on $\Z$: $ X_1,X_2, \ldots ,X_n$, starting from
$l_1< l_2< \ldots <l_n$, respectively.

We observe
these paths between time zero and time $T \in \N$. We define
$\eta (X_1, \ldots, X_j)$ to be the number of walkers
remaining at time $T$ from the initial $j$ at time zero
(i.e., the number of distinct values taken by $X_1(T)$,...,
$X_j(T)$) and $\eta_T = \eta (X_1, \ldots, X_n)$.
We write $X_j-X_i>0$ to denote that $X_j(s)-X_i(s)>0$ for
every $s=0,\ldots,T$.

We are interested in $\P(\eta_T \geq k+1  |\eta_T \geq k)$.
For $k \leq M \leq n$, consider
\[\P(X_1=\xi_1, \ldots, X_{M-1}=\xi_{M-1},Y_M-\xi_{M-1}>0,Y_n-Y_M>0),\]
where $\xi_1, \xi_2, \ldots,\xi_{M-1}$ are  non-crossing
random walk paths starting at $l_1$, $l_2$,$\ldots$, $l_{M-1}$
and $Y_M$ and $Y_n$ are independent simple
random walks starting at $l_M$ and $l_n$ (with no coalescing properties).
With $\xi_1,\xi_2, \ldots,\xi_{M-1}$ fixed
and such that $\eta(\xi_1,\xi_2, \ldots,\xi_{M-1})=k-1$,
we regard only $Y_M$ and $Y_n$ as random and use
the facts that $Y_M$ and $-Y_n$ are independent
and that the events $A = \{Y_M-\xi_{M-1}>0\}$ and
$B=\{Y_n-Y_M>0\}$ are respectively increasing and decreasing
in their dependence on $\{Y_M,-Y_n\}$. It then follows
by the FKG inequalities (for independent variables)
that  $\P(A\cap B)\leq \P(A) \P(B)$.
Letting $B'$ denote the event that $Y_n-Y_1>0$, we obviously
have that $\P(B) \leq \P(B')$ so that
$\P(A\cup B)\leq \P(A)\P(B')$.
Now averaging the latter inequality
over $M, \xi_1,\xi_2, \ldots, \xi_{k-2}$, and noting that
$\P(B')=\P(\eta_T\geq 2)$, we get
$$\P(\eta_T \geq k+1 ) \leq \P(\eta_T \geq k) \P(\eta_T \geq 2)). $$
This completes the proof of the first inequality of
Proposition~\ref{prop:eta}. The second inequality follows immediately and
we leave the equality as an exercise for the reader.

\brm
\label{rm:decay}
The argument in the above proof relies on a negative dependence
(anti-FKG) property
of the point process for which $\hat\eta_{{\cal D}}$ is the cardinality.
This property suggests that the tail decay of $\hat\eta_{{\cal D}}$ is
at at least as fast as Poissonian (which would correspond to the
case of independence).
By analogy with the number of crossings in the scaling limit of
percolation and other statistical mechanics models~\cite{kn:AB},
one could expect the actual decay to be Gaussian.
\erm

The next proposition is a consequence of the one just before; it will
be used in this section and also later in Section~\ref{sec:double}.

\bprop
\label{prop:skel}
Almost surely, for every $\epsilon > 0 $ and every
$\theta = (f(s),t_0)$ in $\bar W({\cal D})$,
there exists a path $
\theta_\epsilon  = (g(s),t'_0)$, in the skeleton
 ${\cal W}({\cal D})$ such that
$g(s) = f(s)$ for all $s \geq t_0+\epsilon$.
\eprop

\noindent{\bf Proof }
Let $\epsilon >0$ be given.
Since $\bw$ is the closure (in $(\Pi,d)$) of
${\cal W}({\cal D})$, we have that there
exists a sequence $\{\theta_n  = (g_n(s),t'_n)\}$ of paths in
${\cal W}({\cal D})$ with $t'_n < t_0 + \frac{\epsilon}{2}$, such that
$d(\theta,\theta_n) \to 0$ as $n \to \infty$. There also exists an integer-valued
random variable $N$ such that $f(t_0 + \frac{\epsilon}{2}) \in [N,N+1]$.
Now by Proposition~\ref{prop:eta} it follows that
$\P(\cup_{k\in \Z} \{ \hat\eta_D (t_0+\frac{\epsilon}{2},\frac{\epsilon}{2}
;k,k+3) < \infty\}) = 1$. Therefore we have
$\hat\eta_D (t_0+\frac{\epsilon}{2},\frac{\epsilon}{2};N-1,N+2)$ is
almost surely
finite. $\theta_n \to \theta$ implies that $ g_n (t_0+\frac{\epsilon}{2})$
is
eventually in $[N-1,N+2]$. Now since the paths are coalescing and $
\hat\eta_D (t_0+\frac{\epsilon}{2},\frac{\epsilon}{2};N-1,N+2)$ is
almost surely
finite we easily see that $\theta_n(t_0+\epsilon)$ is eventually constant.
This implies that $\theta (t_0+\epsilon) = \theta_n(t_0+\epsilon)$ for
large enough $n$ proving the proposition.

\bprop
\label{prop:eta1}
Let $\hat\eta=\hat\eta(t_0,t;a,b)$ be the
counting  random variable for
$\bar \W({\cal D})$, then $\P(\hat\eta \geq k) \leq
{(\Theta(b-a,t))}^k$, and thus $\hat\eta$ is almost surely finite
with finite expectation. Furthermore,
$\hat\eta=\hat\eta_{\cal D}$ almost surely and thus

\beqn
\label{etabound}
\P(\hat\eta\geq k)
&\leq&\P(\hat\eta\geq k-1)\P(\hat\eta\geq 1)\\
&\leq&(\P(\hat\eta\geq 1))^k=(\Theta(b-a,t))^k.
\eeqn

\eprop

\noindent{\bf Proof }
 Choose the set ${\cal D}$ so that its first two points are $(a,t_0)$ and
$(b,t_0)$. We first prove that $\P(\hat\eta \geq k) \leq
{(\Theta(b-a,t))}^k$. From Proposition ~\ref{prop:eta} we know that this claim is
true on the skeleton. Define
$\tilde{\eta}(t_0,t;a,b,\epsilon_1,\epsilon_2)(K)$ = cardinality $\{ y \in
\R :$ there exists a path $(f(s),s_0)$ in $K$ with $s_0 < t_0+\epsilon_2$,
satisfying $f(t_0+\epsilon_2) \in (a-\epsilon_1,b+\epsilon_1)$ and
$f(t_0+t) =y \}$. $\tilde{\eta}_D$ is defined similarly by restricting
the path to be a skeletal path. We observe that $\{K\in \h:\tilde{\eta}_D
\geq k\}$ is an open subset of $\h$. Therefore,
\beqn
\P(\eta(t_0,t;a,b) \geq k) &\leq &\P(\tilde{\eta}
(t_0,t;a,b,\epsilon_1,\epsilon_2) \geq k) \\
\label{eq:eq1}
& = & \P(\tilde{\eta}_D(t_0,t;a,b,\epsilon_1,\epsilon_2) \geq k)\\
\label{eq:eq2}
& \leq & \P(\eta_D(t_0+\epsilon_2,t-\epsilon_2;a-\epsilon_1,b+\epsilon_1)
\geq k) \\
\label{eq:eq3}
& \leq & {(\Theta(b-a+2 \epsilon_1,t+2\epsilon_2))}^{k-1},
\eeqn
 where the equality follows from Proposition~\ref{prop:skel} and
the last inequality follows from inequality~(\ref{eta_D}).
Letting $\epsilon_1 \hbox{ and } \epsilon_2 \to 0$ we see that
\[
 \P(\eta(t_0,t;a,b) \geq k) \leq {(\Theta(b-a,t))}^{k-1}.
\]
Therefore,
\beqn
\P(\hat\eta(t_0,t;a,b) \geq k)  = \P(\eta(t_0,t;a,b) \geq k+1)
\leq  {(\Theta ((b-a),t))}^k.
\eeqn
Using the above inequality for $k = 2$ we get
\beqn
\P(\hat\eta(t_0,t;a,b) \geq 2) \leq{(\Theta ((b-a),t))}^2\leq C {(b-a)}^2, \label{triple}
\eeqn
where $C$ depends only on $t$.
 From inequality~(\ref{triple}) it readily follows that at any
deterministic time almost surely
there are no triple (or higher multiplicity) points. That is,
\[
\P(K\in \h: \hat\eta(t_0,t;c,c) \geq 2 \hbox{ for some c } \in \R)
=0.
\]

Now we want to show that $\hat\eta = \hat\eta_D$ almost surely. If
$\hat\eta \ne \hat\eta_D$, then there exists an $x \in (a,b)$ such that
there is a path $(f(s),t_0) \in \bar\W({\cal D})$ with $f(t_0) = x$ satisfying
the property that if $(g(s),s_0) \in {\cal W}({\cal D})$ with
$s_0 \le t_0$ and
$g(t_0) \in (a,b)$ then $g(s)\ne f(s)$ for all $s\in[t_0,t_0+t]$.

Since ${\cal D}$ is dense in $\R^2$, for all $n \in \N$ we find
paths $(g_n(s), s_n), (h_n(u), u_n)$ in ${\cal W}({\cal D})$
such that $s_n \in (t_0-1/n,t_0), g_n(s_n) \in (x-2/n,x-1/n),
u_n \in (t_0 - 1/n,t_0), h_n(u_n) \in (x+1/n,x+2/n)$. Moreover for
large enough $n$ (depending on $x$) almost
surely we can
choose these paths such that $g_n(t_0) \in (a,x) \hbox{ and }
h_n(t_0) \in (x,b)$. Letting $n \to \infty$  it is not hard to
show that $x$ is
a triple (or higher multiplicity) point.
 To see this, note that by Prop.~\ref{prop:skel}, not only do
$g_n(s)$ and $h_n(s)$ differ from $f(s)$ for
$s \in [t_0,t_0 +t]$ but also for each $s \in (t_0, t_0 +t]$
there are only finitely many possible values for $g_n$ and $h_n$
(as $n$ varies) and hence any (subsequence) limits of
$g_n$ or $h_n$ differ from $f$ for all
$s \in (t_0 +t]$.
Therefore we conclude that
$\P(\hat\eta = \hat\eta_D) = 1$, thus proving the proposition.

\bteo
\label{teo:char1}
Let ${\cal W}'$ be an $({\cal H},{\cal F}_{{\cal H}})$-valued random
variable; its distribution
equals that of the (standard) Brownian web $\bw$ (as characterized
by Theorems \ref{teo:char} and \ref{teo:charm})
if its finite dimensional distributions are coalescing (standard)
Brownian motions (i.e.
conditions (o) and (i) of Theorem~\ref{teo:char} are valid) and
\begin{itemize}
\item[$(ii'')$] for all $t_0,t,a,b,\,\,\hat\eta_{{\cal W}'}$ is equidistributed
with $\hat\eta_{\bw}$.
\end{itemize}
\eteo

For purposes of proving our convergence results, we will use a
modified version of the above
characterization theorem in which conditions $(o)$, $(i)$, $(ii'')$ are all
weakened.

\bteo
\label{teo:char2}
Let ${\cal W}'$ be an $({\cal H},{\cal F}_{{\cal H}})$-valued random
variable and let $\cal D$ be
a countable dense deterministic subset of $\r^2$ and for each
$y\in\cal D$, let $\theta^y\in {\cal W}'$
be some single (random) path starting at $y$. ${\cal W}'$ is equidistributed
with the (standard) Brownian
web $\bar {\cal W}$ if
\begin{itemize}
\item[(i$^\prime$)] the $\theta^y$'s are distributed as coalescing
(standard)
Brownian motions and
\item[(ii$^{\prime\prime\prime}$)] for all $t_0,t,a,b,\,\,
\eta_{{\cal W}'}<<\eta_{\bar{\cal W}}$,
i.e.~$\P(\eta_{{\cal W}'}\geq
k)\leq\P(\eta_{\bar {\cal W}}\geq k)$ for all $k$.
\end{itemize}
\eteo

\noindent {\bf Proof } We need to show that the above conditions
together imply that $\mu'$, the distribution of $\W'$, equals the
distribution $\mu$ of the constructed Brownian web $\bar{\cal W}$. Let
$\eta'$ the counting random variable appearing in condition
(ii$^{\prime\prime\prime}$) for $\mu'$.
Choose some deterministic dense countable subset $\cd$ and consider
the countable
collection ${\cal W}^*$ of paths of ${\cal W}'$ starting from $\cd$. By condition
(i$^\prime$), ${\cal W}^*$
is equidistributed with our constructed
Brownian skeleton ${\cal W}$ (based on the same $\cd$)
and hence the closure $\bar{{\cal W}^*}$ of ${\cal W}^*$ in $(\o,d)$ is a subset of
${\cal W}'$
that is equidistributed with our constructed Brownian
web $\bar{\cal W}$. To complete the proof, we will use condition
(ii$^{\prime\prime\prime}$) to show
that ${\cal W}'\setminus \bar{{\cal W}^*}$ is almost surely empty
by using the fact that the
counting  random variable $\eta^*$ for $\bar{{\cal W}^*}$
already satisfies
condition (ii$^{\prime\prime\prime}$) since $\bar{{\cal W}^*}$ is
distributed as a Brownian web. If
${\cal W}'\setminus \bar{\cal W}^*$ were nonempty (with strictly positive
probability),
then there would have to be some rational $t_0, t, a, b$
for which $ \eta' > \eta^*$. But then
\begin{equation}
\P(\eta'(t_0,t;a,b) > \eta^*(t_0,t;a,b)) \,>\, 0
\end{equation}
for some rational $t_0, t, a, b$, and
this together with the fact that $\P(\eta'\geq\eta^\ast)=1$ (which follows from
$\bar\W^\ast\subset\W'$) would violate condition (ii$^{\prime\prime\prime}$)
with those  $t_0, t, a, b$. The proof is complete.

\brm
The condition $\eta_{{\cal W}'}<<\eta_{\bar{\cal W}}$ can be replaced by
$\E(\eta_{{\cal W}'})\leq\E(\eta_{\bar{\cal W}})$.
We note that $\E(\eta_{\bar{\cal W}})=1+(b-a)/\sqrt{\pi t}$, as given
in~\cite{kn:FINR} by a calculation stretching back to~\cite{kn:BG}.
 So, in particular, in the context of Theorem~\ref{teo:char2}, if
besides ($i'$),
$\E(\eta_{{\cal W}'})\leq1+(b-a)/\sqrt{\pi t}$ for all $t_0,t,a,b$, then
${\cal W}'$ is equidistributed with the Brownian web.
\erm

%%%%%%%%%%%%%%%%%%%%%%%%%%%%%%%%%%%%%%%%%%%%
%%%%%%%%%%%%%%% SECTION 5 %%%%%%%%%%%%%%%%%%%%%%%
%%%%%%%%%%%%%%%%%%%%%%%%%%%%%%%%%%%%%%%%%%%%
%%%%%%%%%%%%%%%%%%%%%%%%%%%%%%%%%%%%%%%%%%%%

\section{Dual and Double Brownian Webs}

\setcounter{equation}{0}
\label{sec:double}

\begin{figure}[!ht]
\begin{center}
\includegraphics[width=10cm]{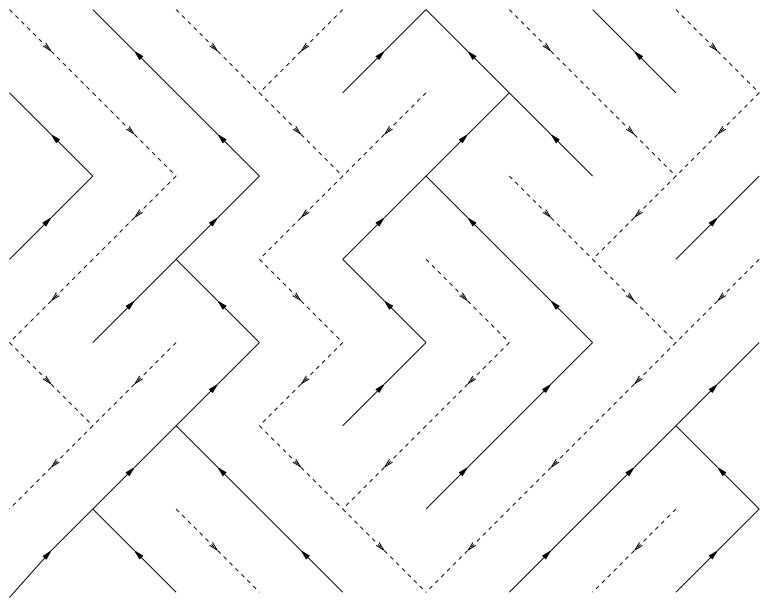}
\caption{Forward coalescing random walks (full lines) in discrete time and their
dual backward walks (dashed lines).}\label{figdouble}
\label{double}
\end{center}
\end{figure}

In this section, we construct and characterize the {\em Double
Brownian Web}, which combines the
Brownian web with a {\em Dual Brownian Web} of coalescing Brownian
motions moving backwards in time.
\brm
\label{rmk:dual}
In the graphical representation of Harris
for the one-di\-men\-sio\-nal voter model~\cite{kn:H},
coalescing random walks forward in time and coalescing dual random walks
backward in time (with forward and backward walks not crossing each
other)
are constructed simultaneously (see, e.g., the discussion
in~\cite{kn:FINS1,kn:FINS2}).
Figure~\ref{figdouble} provides an example in discrete time.
Note that there is no crossing between forward and backward walks---a
property that is preserved in the
{\em Double Brownian Web} (DBW) scaling limit.
The simultaneous construction of forward and (dual) backward Brownian
motions
was emphasized in~\cite{kn:TW,kn:STW} and their approach and results
can be
applied to extend both our characterization and convergence results to
the DBW which includes simultaneously the
forward BW
and its dual backward BW. We note that in the DBW, the $\eta$ of
Definition~\ref{eta}
equals (almost surely for deterministic $t_0,t,a,b$)
$1+ \eta^{\mbox{\rm\scriptsize dual}}$, where
$\eta^{\mbox{\rm\scriptsize dual}}$ is the number of distinct points in
$[a,b] \times \{t_0\}$ touched by backward paths which also touch
$\R \times \{t_0 + t\}$.
\erm

%%%%%%%%%%%%%%%%%%%%%%%%%%%%%%%%%%%%%%%%%%%%%%
%%%%%%%%%%%%%%  FIGURE 2 %%%%%%%%%%%%%%%%%%%%%%%%%%
%%%%%%%%%%%%%%%%%%%%%%%%%%%%%%%%%%%%%%%%%%%%%%

The forward and backward paths basically reflect off each other. In
Section~\ref{sec:rw-conv} (see Theorem~\ref{teo:convdbl}),
we show that the
Double Brownian Web arises as the scaling limit of coalescing random
walks forward in time together
with their (graph theoretic) dual coalescing paths backward in time
(see Figure~\ref{figdouble}).
Our analysis will rely on a paper~\cite{kn:STW}
of Soucaliuc, T\'oth and Werner together with our results from earlier
sections on the (forward)
Brownian web.

We again begin with an (ordered) dense countable set ${\cal
D}\subset\r^2$ but this time we need
twice as many i.i.d.~standard B.M.'s as before
$B_1,B_1^b,B_2,B_2^b,\ldots$. We use this paths
to construct forward and backward independent B.M.'s
$W_1,W_1^b,W_2,W_2^b,\ldots$ starting from
$(x_j,t_j)\in{\cal D}$ as follows.
\begin{eqnarray}
  \label{eq:for}
  W_j(t)&=&x_j+B_j(t-t_j),\,t\geq t_j\\
  \label{eq:bac}
  W_j^b(t)&=&x_j+B_j^b(t_j-t),\,t\leq t_j.
\end{eqnarray}
Then we create coalescing and reflecting Brownian paths
$\tilde W_1,\tilde W_1^b,\ldots$ inductively, as follows.
 \begin{eqnarray}
  \label{eq:tw1}
  \tilde W_1&=&W_1;\quad\tilde W_1^b\,\,=\,\,W_1^b;\\
  \label{eq:twn}
  \tilde W_n&=&
  CR(W_n;\tilde W_1,\tilde W_1^b,\ldots,\tilde W_{n-1},\tilde W_{n-1}^b);\\
  \label{eq:twnd}
  \tilde W_n^b&=&
  CR(W_n^b;\tilde W_1,\tilde W_1^b,\ldots,\tilde W_{n-1},\tilde W_{n-1}^b),
\end{eqnarray}
where the operation $CR$ is defined in~\cite{kn:STW}, Subsubsection 3.1.4.
We proceed to explain $CR$ for the simplest case, in the definition of $\tilde W_2$.

As pointed out in~\cite{kn:STW}, the nature of the reflection of a
forward Brownian path $\tilde W$ off a backward Brownian path $\tilde W^b$
(or vice-versa) is
special. It is actually better described as a push of $\tilde W$ off $\tilde W^b$
(see Subsection 2.1 in~\cite{kn:STW}). It does not have an explicit
form in general, but in the case of one forward path and one backward
path, the form is as follows. Following our notation and construction,
we ignore $\tilde W_1$ and consider $\tilde W_1^b$ and  $\tilde W_2$
in the time interval $[t_2,t_1]$ (we suppose $t_2<t_1$; otherwise,
$\tilde W_1^b$ and  $\tilde W_2$ are independent). Given $W_2$ and $\tilde W_1^b$,
for $t_2\leq t\leq t_1$,
\begin{equation}
  \label{eq:refl}
  \tilde W_2(t)=
  \begin{cases}
    W_2(t)+\sup_{t_2\leq s\leq t}(W_2(s)-\tilde W_1^b(s))^-,
          \mbox{ if } W_2(t_2)>\tilde W_1^b(t_2);\\
    W_2(t)-\sup_{t_2\leq s\leq t}(W_2(s)-\tilde W_1^b(s))^+,
          \mbox{ if } W_2(t_2)<\tilde W_1^b(t_2).
  \end{cases}
\end{equation}
After $t_1$, $\tilde W_2$ interacts only with $\tilde W_1$,
by the usual coalescence.

We call
${\cal W}^D_n:=\{\tilde W_1,\tilde W_1^b,\ldots,\tilde W_n,\tilde W_n^b\}$
{\em coalescing/reflecting forward and backward Brownian motions (starting at
$\{(x_1,t_1),\ldots,(x_n,t_n)\}$)}.
We will also use the alternative notation ${\cal W}^D({\cal D}_n)$ in place of
${\cal W}^D_n$, where ${\cal D}_n:=\{(x_1,t_1),\ldots,(x_n,t_n)\}$.
\brm
\label{rm:stw}
In Theorem 8 of~\cite{kn:STW}, it is proved that the above
construction is a.s.~well-defined, gives a perfectly
coalescing/reflecting system (see Subsubsection 3.1.1 in~\cite{kn:STW}),
and for every $n\geq1$, the distribution of ${\cal W}^D_n$
does not depend on the ordering of ${\cal D}_n$.
It also follows from that result that $\{\tilde W_1,\ldots,\tilde W_n\}$ and
$\{\tilde W_1^b,\ldots,\tilde W_n^b\}$ are separately forward and backward coalescing
Brownian motions, respectively. Thus
$\{\tilde W_1,\tilde W_2,\ldots\}$ and
$\{\tilde W_1^b,\tilde W_2^b,\ldots\}$ are forward and backward
Brownian web skeletons, respectively.
\erm

\brm
One can alternatively use a set ${\cal D}^b$ of starting points for the backward
paths different than ${\cal D}$
rather than our choice above of ${\cal D}^b={\cal D}$.
\erm

We now define dual spaces of paths going backward in time $(\Pi^b,d^b)$ and
a corresponding $({\cal H}^b,d_{{\cal H}^b})$
in an obvious way, so that they are the dual versions of
$(\Pi,d)$ and $({\cal H},d_{{\cal H}})$,
and then define ${\cal H}^D={\cal H}\times{\cal H}^b$ and
$$d_{{\cal H}^D}((K_1,K_1^b),(K_2,K_2^b))=\max(d_{{\cal
H}}(K_1,K_2),d_{{\cal H}^b}(K_1^b,K_2^b)).$$

As in the construction of the (forward) BW, we now define
\beqn
{\cal W}_n^D({\cal D})
&=&
\{\tilde W_1,\ldots,\tilde W_n\}\times\{\tilde W_1^b,\ldots,\tilde W_n^b\},\\
{\cal W}^D({\cal D})
&=&\{\tilde W_1,\tilde W_2,\ldots\}\times\{\tilde W_1^b,\tilde W_2^b,\ldots\},\\
\bar{\cal W}^D({\cal D})&=&
\overline{\{\tilde W_1,\tilde W_2,\ldots\}}\times\overline{\{\tilde W_1^b,\tilde W_2^b,\ldots\}}.
\eeqn
The latter closures are in $\Pi$ for the first factor and in $\Pi^b$
for the second one.

{}From Remark~\ref{rm:stw}, we have that
$$\bar{\cal W}:=\overline{\{\tilde W_1,\tilde W_2,\ldots\}} \mbox{ and }
\bar{\cal W}^b:=\overline{\{\tilde W_1^b,\tilde W_2^b,\ldots\}}$$
are forward and backward Brownian webs, respectively. The next
result follows from Proposition~\ref{prop:bwcpt}.
\bprop
\label{prop:dbwcpt}
Almost surely, $\bar{\cal W}^D({\cal D})\in{\cal H}^D$  (i.e.
$\overline{\{\tilde W_1,\tilde W_2,\ldots\}}$ and
$\overline{\{\tilde W_1^b,\tilde W_2^b,\ldots\}}$ are compact).
\eprop

\brm
\label{rm:dlim} It is immediate from this proposition that
$$\bar{\cal W}^D({\cal D})=\lim_{n\to\infty}{\cal W}_n^D({\cal D}),$$
\erm
where the limit is in the $d_{{\cal H}^D}$ metric.

\bprop
\label{prop:dpaths}
$\bar{\cal W}^D({\cal D})$ satisfies
\begin{itemize}
\item[$(o^D)$] From any deterministic $(x,t)$ there is almost surely
a unique forward path and unique backward path.
\item[$(i^D)$] For any deterministic
${\cal D}'_n:=\{(y_1,s_1),\ldots,(y_n,s_n)\}$ the
forward and backward paths
from ${\cal D}'_n$, denoted $\bar{\cal W}^D({\cal D},{\cal D}'_n)$,
are distributed as coalescing/reflecting forward and backward Brownian
motions starting at ${\cal D}'_n$. In other words,
$\bar{\cal W}^D({\cal D},{\cal D}'_n)$ has the same distribution as
${\cal W}^D({\cal D}'_n)$.
\end{itemize}
\eprop

\noindent{\bf Proof } $(o^D)$ follows from the separate properties
of $\bar{\cal W}({\cal D})$ and of $\bar{\cal W}^b({\cal D})$,
since they are distributed respectively as forward and backward
Brownian webs, as previously noted.

As to $(i^D)$: for $k=1,2,\ldots$ let
${\cal D}^{(k)}_n=\{(x^{(k)}_1,t^{(k)}_1),\ldots,(x^{(k)}_n,t^{(k)}_n)\}
\subset{\cal D}$
be such that $(x^{(k)}_i,t^{(k)}_i)$ converges to
$(y_i,s_i)$ for each $1\leq i\leq n$ as $k\to\infty$.
Since ${\cal W}^D({\cal D}'_n)$ and ${\cal W}^D({\cal D}^{(k)}_n)$
(for every $k\geq1$) are defined using the same $(B_i,B_i^b)$'s,
we get from standard Brownian motion sample path properties that
${\cal W}^D({\cal D}^{(k)}_n)\to{\cal W}^D({\cal D}'_n)$ a.s.~as
$k\to\infty$, and thus
\begin{equation}
  \label{eq:dp1}
 {\cal W}^D({\cal D}^{(k)}_n)\to{\cal W}^D({\cal D}'_n)
\end{equation}
in distribution as $k\to\infty$. Now from Theorem 8 of~\cite{kn:STW},
$\bar{\cal W}^D({\cal D},{\cal D}^{(k)}_n)$ has the same distribution as
$\bar{\cal W}^D({\cal D}^{(k)}_n)$. We proceed to show that
\begin{equation}
  \label{eq:dp2}
  \bar{\cal W}^D({\cal D},{\cal D}^{(k)}_n)\to
  \bar{\cal W}^D({\cal D},{\cal D}'_n)
\end{equation}
a.s.~as $k\to\infty$, and $(i^D)$ follows immediately from that
and~(\ref{eq:dp1}).

Denoting by $\bar{\cal W}({\cal D},{\cal D}^{(k)}_n)$
and $\bar{\cal W}^b({\cal D},{\cal D}^{(k)}_n)$ respectively the set
of forward and backward paths of
$\bar{\cal W}^D({\cal D},{\cal D}^{(k)}_n)$,
(\ref{eq:dp2}) is
equivalent to the pair of a.s.~limits (as $k\to\infty$):
\begin{eqnarray}
  \label{eq:dp3a}
  \bar{\cal W}({\cal D},{\cal D}^{(k)}_n)&\to&
  \bar{\cal W}({\cal D},{\cal D}'_n);\\
  \label{eq:dp3b}
  \bar{\cal W}^b({\cal D},{\cal D}^{(k)}_n)&\to&
  \bar{\cal W}^b({\cal D},{\cal D}'_n).
\end{eqnarray}
We argue~(\ref{eq:dp3a}) only;
(\ref{eq:dp3b}) is similar. Let $\g^{(k)}_i$ be the path in
$\bar{\cal W}({\cal D},{\cal D}^{(k)}_n)$ starting at
$(x^{(k)}_i,t^{(k)}_i)$ and $\g_i$ be the path in
$\bar{\cal W}({\cal D},{\cal D}'_n)$ starting at
$(y_i,s_i)$. It is enough to show that for every
$i=1,\ldots,n$, $\g^{(k)}_i\to\g_i$ in $\Pi$ as $k\to\infty$.
Suppose not; then since $\g^{(k)}_i\in\bar{\cal W}({\cal D})$, which is compact,
there must be a subsequence of $(\g^{(k)}_i)_k$ which converges to some
$\g'\ne\g_i$. Since $\g'\in\bar{\cal W}({\cal D})$,
this would contradict $(o^D)$. The proof is finished.

\bprop
\label{prop:dind}
The distribution of $\bar{\cal W}^D({\cal D})$ as an
$({\cal H}^D,{\cal F}_{{\cal H}^D})$-valued
random variable (where
${\cal F}_{{\cal H}^D}={\cal F}_{{\cal H}}\times{\cal F}_{{\cal H}^b}$),
does not depend on ${\cal D}$.
Furthermore,
\begin{itemize}
\item[$({ii}^{D})$] for any deterministic dense ${\cal D}'$, almost surely
\begin{equation*}
\bar{\cal W}^D({\cal D})=
\overline{\{W_{x,t}:(x,t)\in{\cal
D}'\}}\times\overline{\{W_{x,t}^b:(x,t)\in{\cal D}'\}},
\end{equation*}
\end{itemize}
where $W_{x,t},W_{x,t}^b$ are respectively the forward and backward
paths in $\bar{\cal W}^D({\cal D})$
starting from $(x,t)$, and
the closures in $({ii}^{D})$ are in $\Pi$ for the first factor and in $\Pi^b$ for the second one.
\eprop

\noindent{\bf Proof } The proof is essentially the same
as that of Proposition~\ref{prop:ind}, using
Proposition~\ref{prop:dpaths} in place of Proposition~\ref{prop:findim}
and the natural (and similarly proved) analogue of
Proposition~\ref{prop: D_1-D_2}.

We now give some characterization theorems for the
distribution of the double Brownian web
$\bar\W^D$.

\bteo
\label{teo:chard}
The double Brownian web is characterized
(in distribution, on $({\cal H}^D, {\cal F}_{{\cal H}^D})$) by conditions
$(o^D),(i^D)$ and $(ii^D)$.
\eteo

\noindent{\bf Proof }
The proof is basically the same as for Theorem~\ref{teo:char}.
We include it for completeness.
Take any
\(({\cal H}^D, {\cal F}_{{\cal H}^D}) \)-valued random variable $\X^D$ with properties
$(o^D), (i^D)$ and $(ii^D)$, and fix a deterministic countable dense subset
${\cal D}=\{(x_1,t_1),(x_2,t_2),\ldots\}$ of $\r^2$.
Let $\bar\W^D=\bar\W^D({\cal D})$ be the version of the double
Brownian web constructed
earlier in this section of the paper, using $\cal D$
(and recall that the distribution
does not depend on $\cal D$, by Proposition~\ref{prop:dind}).
 From $(i^D)$ for $\X^D$ and Proposition~\ref{prop:dpaths},
$\X_n^D:=\{X_{x_i,t_i},X^b_{x_i,t_i};\,i=1,\ldots,n\}$, where
$X_{x,t},X_{x,t}^b$ are the forward and backward paths of $\X^D$ starting at $(x,t)$, respectively
(almost surely a unique pair by $(o^D)$),
is equidistributed with
$\W_n^D=\{W_{x_i,t_i},W^b_{x_i,t_i};\,i=1,\ldots,n\}$,
where
$W_{x,t},W_{x,t}^b$ are the forward and backward paths of $\bar\W^D$ starting at $(x,t)$,
for all $n\geq1$.
Now $(ii^D)$ for $\X^D$
and Remark~\ref{rm:dlim} imply or say that $\X^D=\lim_{n\to\infty}\X^D_n$ and
$\bar\W^D=\lim_{n\to\infty}\W^D_n$ a.s.~.
So, they have the same distribution and the proof is complete.

\bdf
(Stochastic ordering)
We define a partial order in ${\cal H}^D$ as follows:
given $(K_1,K_1^b),\,(K_2,K_2^b)\in{\cal H}^D$, we say that
$(K_1,K_1^b)\leq(K_2,K_2^b)$ if $K_1\subseteq K_2$ and $K_1^b\subseteq
K_2^b$.
For two measures $\mu_1^D,\,\mu_2^D$ in $({\cal H}^D,{\cal F}_{{\cal
H}^D})$,
we say that $\mu_1^D<<\mu_2^D$ if $\int f\,d\mu_1^D\leq\int f\,d\mu_2^D$
for every bounded $f\in{\cal F}_{{\cal H}^D}$ that is increasing in the
partial
order $\leq$.
\edf

\bteo
\label{teo:charmd}
The distribution $\mu_{\bar\W^D}$ of the
double Brownian web $\bar\W^D$ is characterized by conditions
$(o^D),(i^D)$ and
\begin{itemize}
\item[$({ii'}^{D})$] any $({\cal H}^D,{\cal F}_{{\cal
H}^D})$-valued random variable, $\X^D$,
satisfying $(o^D)$ and $(i^D)$ has $\mu_{\bar\W^D}<<\mu_{\X^D}$.
\end{itemize}
\eteo

\noindent{\bf Proof }
The proof is analogous to the one of Theorem~\ref{teo:charm}.
Suppose ${\bar\W^D}$ is an $({\cal H}^D,{\cal F}_{{\cal H}^D})$-valued
random variable satisfying $(o^D)$, $(i^D)$ and $({ii'}^{D})$.
Fix a deterministic countable dense subset
${\cal D}=\{(x_1,t_1),(x_2,t_2),\ldots\}$ of $\r^2$ and let $\X^D$ be the double Brownian web
constructed from $\cal D$ (denoted earlier by $\bar\W^D(\cal D))$.
Since $\bar\W^D\supset\overline{\{W_{x,t},W^b_{x,t};\,(x,t)\in{\cal D}\}}$,
where $W_{x,t},W^b_{x,t}$ are the forward and backward paths of $\bar\W^D$ starting at $(x,t)$,
and, from $(o^D)$ and $(i^D)$ for ${\bar\W^D}$,
$\W_n^D=\{W_{x_i,t_i},W^b_{x_i,t_i};\,i=1,\ldots,n\}$
has the same distribution
as $\X_n^D:=\{X_{x_i,t_i},X^b_{x_i,t_i};\,i=1,\ldots,n\}$ for all $n\geq1$, where
$X_{x,t},X^b_{x,t}$ are the forward and backward paths of $\X^D$ starting at $(x,t)$,
we have that
\( \mu_{{\X}^D} << \mu_{{\bar\W}^D} \). This and $(ii'^D)$ together imply that
\( \mu_{{\X}^D} = \mu_{{\bar\W}^D} \), and the proof is complete.

\bdf
\label{etab}
For $t>0,\,t_0,a,b\in\r,\,a<b$, let $\eta^b(t_0,t;a,b)$ be the number
of {\em distinct} points in
$\R\times\{t_0-t\}$ that are touched by paths in $\bar W^b$ which also
touch
some point in $[a,b]\times\{t_0\}$. Let also
${\hat\eta}^b(t_{0}, t; a,b)=\eta^b(t_{0}, t; a,b)-1.$
\edf

\bteo
\label{teo:chard2}
Let ${\W'}^{D}=(\W',{\W'}^{b})$ be an $({\cal H}^D,{\cal F}_{{\cal
H}^{D}})$-valued random variable,
let ${\cal D}$ be
a countable dense deterministic subset of $\r^2$, and for each
$y\in{\cal D}$ let
$\theta^y\in \W'$ and ${\theta^y}^b\in {\W'}^b$ be single paths starting
at $y$.
${\W'}^D$ is equidistributed with the double Brownian web
$\W^D = (\bar \W,\bar \W^b)$ if
\begin{itemize}
\item[$({i'}^D)$]
the $\theta^y$'s and ${\theta^y}^b$'s are distributed as
coalescing/reflecting
forward and backward Brownian motions, and
\item[$({ii'}^D)$] for all $t_0,t,a,b,
\,\, \hat\eta_{\W'}<<\hat\eta_{\bar \W}$ and
$\hat\eta_{{\W'}^b}<<\hat\eta_{\bar \W^b}$.
\end{itemize}
\eteo

\noindent{\bf Proof }
The proof is analogous to the one of Theorem~\ref{teo:char2}.
We need to show that the above conditions
together imply that $\mu'$, the distribution of ${\W'}^D$, equals the
distribution $\mu$ of the previously
constructed double Brownian web $\bar{\cal W}^D$. Let
$\eta',{\eta'}^b$ denote the forward and backward counting random variables
appearing in condition
(ii$^{\prime\prime\prime}$) for $\mu'$.
Choose some deterministic dense countable subset $\cd$ and consider
the countable
collection ${{\cal W}_*}^D$ of forward and backward paths of ${{\cal W}'}^D$
starting from $\cd$. By $({i'}^D)$, ${{\cal W}_*}^D$
is equidistributed with our constructed
double Brownian web skeleton ${\cal W}^D$ (based on the same $\cd$)
and hence the closure ${\bar{\cal W}_*}^D$ of ${{\cal W}_*}^D$  is a subset of
${{\cal W}'}^D$
that is equidistributed with our constructed double Brownian
web $\bar{\cal W}^D$. To complete the proof, we will use condition
$({ii'}^D)$ to show
that ${{\cal W}'}^D\setminus {\bar{\cal W}_*}^D$ is almost surely empty,
by using the fact that
the forward and backward
counting  random variables $\eta_*,{\eta_*}^b$ for ${\bar{\cal W}_*}^D$
already satisfy
condition $({ii'}^D)$, since ${\bar{\cal W}_*}^D$ is
distributed as a double Brownian web. If
${{\cal W}'}^D\setminus {\bar{\cal W}_*}^D$ were nonempty (with strictly positive
probability),
then there would have to be some rational $t_0, t, a, b$
for which either $ \eta' > \eta_*$ or $ {\eta'}^b > {\eta_*}^b$. But then
\begin{equation}
\P\left(\eta'(t_0,t;a,b)>\eta_*(t_0,t;a,b)\mbox{ or }{\eta'}^b(t_0,t;a,b)>{\eta_*}^b(t_0,t;a,b)\right)>0
\end{equation}
for some rational $t_0, t, a, b$, and
this together with the fact that $\P(\eta'\geq\eta_\ast,\,{\eta'}^b\geq{\eta_\ast}^b)=1$ (which follows from
${\bar\W_\ast}^D\subset{\W'}^D$) would violate condition $({ii'}^D)$
with those  $t_0, t, a, b$. The proof is complete.

We now discuss ``types'' of points $(x,t)\in\r^2$,
whether deterministic or not.
For the (forward) Brownian web,
we define
\beqn\nn
m_{\mbox{{\scriptsize in}}}(x_0,t_0)\!\!&=&\!\!
\lim_{\epsilon\downarrow0}\{\mbox{number of paths in } \W
\mbox{ starting at some }
t_0-\epsilon
\mbox{ that pass}\\
\label{eq:min}
&& \mbox{through } (x_0,t_0) \mbox{ and are disjoint for }
t_0-\epsilon<t<t_0\};\\\nn
m_{\mbox{{\scriptsize out}}}(x_0,t_0)\!\!&=&\!\!
\lim_{\epsilon\downarrow0}\{\mbox{number of paths in } \W
\mbox{ starting at }
(x_0,t_0)\mbox{ that are}\\
\label{eq:mout}
&& \mbox{ disjoint for }
t_0<t<t_0+\epsilon\}.
\eeqn

For $\W^b$, we similarly define $m^b_{\mbox{{\scriptsize in}}}(x_0,t_0)$ and
$m^b_{\mbox{{\scriptsize out}}}(x_0,t_0)$.

\begin{figure}[!ht]
\begin{center}
\includegraphics[width=6cm]{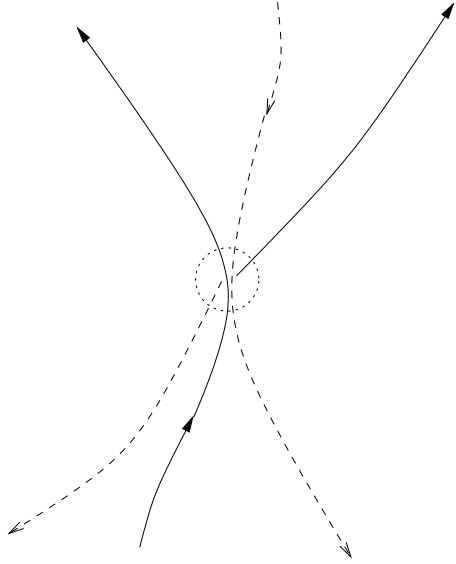}
\caption{A schematic diagram of a point $(x_0,t_0)$ of type
$(m_{\mbox{{\scriptsize\emph{in}}}},m_{\mbox{{\scriptsize \emph{out}}}})=(1,2)$,
with necessarily also
$(m^b_{\mbox{{\scriptsize\emph{in}}}},m^b_{\mbox{{\scriptsize \emph{out}}}})=(1,2)$.
In this example the incoming forward path connects to the leftmost outgoing path
(with a corresponding dual connectivity for the backward paths);
at some of the other points of type $(1,2)$ it
will connect to the rightmost path.}\label{fig12}
\label{12} 
\end{center}
\end{figure}

\bdf
The type of $(x_0,t_0)$ is the pair $(m_{\mbox{{\scriptsize
\emph{in}}}},
m_{\mbox{{\scriptsize \emph{out}}}})$---see Figure~\ref{fig12}.
We denote by $S_{i,j}$ the set of points of $\r^2$ that are of type
$(i,j)$, and by $\bar S_{i,j}$ the set of points of $\r^2$ that are of
type $(k,l)$ with $k\geq i$, $l\geq j$.
\edf

%%%%%%%%%%%%%%%%%%%%%%%%%%%%%%%%%%%%%%%%%%%%%%
%%%%%%%%%%%%%  FIGURE 3 %%%%%%%%%%%%%%%%%%%%%%%%%%%
%%%%%%%%%%%%%%%%%%%%%%%%%%%%%%%%%%%%%%%%%%%%%%

\brm
\label{rm:dense}
Using the translation and scale invariance properties of the
Brownian web distribution, it can be shown that
for any $i,j$, whenever $S_{i,j}$ is nonempty, it must be dense in $\r^2$.
The same can be
said of $S_{i,j}\cap\r\times\{t\}$ for deterministic $t$.
These denseness properties can also be shown for each $i,j$ by more
direct arguments.
\erm

\bprop
\label{prop:dual}
For the double Brownian web, almost surely for {\em every}
$(x_0,t_0)\in\r^2$, $m^b_{\mbox{{\scriptsize \emph{in}}}}
(x_0,t_0)=m_{\mbox{{\scriptsize \emph{out}}}}(x_0,t_0)-1$
and
$m^b_{\mbox{{\scriptsize \emph{out}}}}(x_0,t_0)=
m_{\mbox{{\scriptsize \emph{in}}}}(x_0,t_0)+1$. See Figure~\ref{fig12}.
\eprop

\noindent{\bf Proof } It is enough to prove
(i) that for every incoming forward path to a point $(x,t)$,
there are two locally disjoint backward paths starting at that point
with one on either side of the forward path;
and (ii) that for every two locally disjoint backward paths
starting at a point $(x,t)$, there is an
incoming forward path to $(x,t)$ between the two backward paths.
(Note that by a $t \longleftrightarrow -t$ time reflection
argument, one would then get a similar result for incoming backward
paths and pairs of outgoing forward paths.)

Let us start with the first assertion. Let $\g$ be an incoming
forward path to $(x,t)$. This means that the
starting time $s$ of $\g$ is such that $s<t$.
By Proposition~\ref{prop:skel}, the portion of $\g$ above time $s+\e$
is in the forward skeleton for every $\e>0$. Now consider a
sequence of pairs of backward paths $(\g_k,\g'_k)$ starting at
$((x_k,t_k),(y_k,s_k))\in{\cal D}\times{\cal D}$ with
$((x_k,t_k),(y_k,s_k))\to((x,t),(x,t))$ as $k\to\infty$,
$s+\e < s_k,t_k<t$, $x_k<\g(t_k)$ and $y_k>\g(s_k)$. From the reflection
of the forward and backward skeletons off each other and
the fact that two backward paths in the skeleton must coalesce
once they meet, it follows that
$\g_k(t')<\g'_k(t')$ for all $t'\in[\mbox{max}\{s_k,t_k\},t']$.
We then conclude from
compactness that there are two locally disjoint limit paths,
one for $(\g_k)$ and one for $(\g'_k)$, both starting from $(x,t)$.

We argue (ii) similarly. Given two locally disjoint backward paths
$\g,\g'$ starting at $(x,t)$, there exists $s<t$ such that
either $\g(t')<\g'(t')$ for $s<t'<t$ or
$\g'(t')<\g(t')$ for $s<t'<t$. Suppose it is the first case;
otherwise, switch labels. Then choose a point $(x',s')\in{\cal D}$
with $s<s'<t$ and $\g(s')<x'<\g'(s')$. The fact that the portions of
$\g$ and $\g'$ below time $t-\e$
is in the backward skeleton for every $\e>0$ and the reflection
of the forward and backward skeletons off each other now implies that
the forward path starting at $(x',s')$ is squeezed between
$\g$ and $\g'$ and goes to $(x,t)$.

\bteo
\label{teo:types}
For the (double) Brownian web, almost surely, every $(x,t)$ has one of the
following types, all of which occur: $(0,1)$, $(0,2)$, $(0,3)$, $(1,1)$,
$(1,2)$, $(2,1)$.
\eteo

\brm
Points of type $(1,2)$ are particularly
interesting in that the single incident path continues along exactly one
of the two outward paths --- with the choice determined intrinsically
rather than by some convention. See Figure~\ref{fig12} for a schematic
diagram of a ``left-handed'' continuation. An $(x_0,t_0)$ is of
type $(1,2)$ precisely if both a forward and a backward path pass through
$(x_0,t_0)$. It is either left-handed or right-handed according to
whether the forward path is to the left or the right of the backward
path near $(x_0,t_0)$. Both varieties occur and the proof of
Theorem~\ref{teo:typesa} below shows that the Hausdorff dimension
of $1$ applies separately to each of the two varieties.
\erm

T\'oth and Werner~\cite{kn:TW} gave a definition of types of points of
$\r^2$ similar to ours, but for a somewhat different process
and proved the above theorem with that definition and for that process
(see definition at page 385, paragraph of equation~(2.28) and
Proposition 2.4 in~\cite{kn:TW}).
One way then to establish Theorem~\ref{teo:types} is to show the equivalence
of ours and T\'oth and Werner's definition and that their arguments
hold for our process. We prefer, for the sake of simplicity and
completeness, to give a direct argument, out of which the following
complementary results also follow.

\bteo
\label{teo:typesa}
Almost surely,
$S_{0,1}$ has full Lebesgue measure in $\r^2$,
$S_{1,1}$ and $S_{0,2}$ have Hausdorff dimension $3/2$ each,
$S_{1,2}$ has Hausdorff dimension $1$, and
$S_{2,1}$ and $S_{0,3}$ are both countable and dense in $\r^2$.
\eteo

\bteo
\label{teo:typesb}
Almost surely: for every $t$
\begin{itemize}
\item[a)] $S_{0,1}\cap\r\times\{t\}$ has full Lebesgue measure in $\r\times\{t\}$;
\item[b)] $S_{1,1}\cap\r\times\{t\}$ and $S_{0,2}\cap\r\times\{t\}$
are both countable and dense in $\r\times\{t\}$;
\item[c)] $S_{1,2}\cap\r\times\{t\}$, $S_{2,1}\cap\r\times\{t\}$ and
$S_{0,3}\cap\r\times\{t\}$ have all cardinality at most $1$.
\end{itemize}

For every deterministic $t$,
$S_{1,2}\cap\r\times\{t\}$, $S_{2,1}\cap\r\times\{t\}$ and
$S_{0,3}\cap\r\times\{t\}$ are almost surely empty.
\eteo

\noindent{\bf Proof of Theorems~\ref{teo:types} and~\ref{teo:typesa} }
We start by ruling out the cases that do not occur almost surely.
For $i,j\geq0$, $S_{i,j}=\emptyset$ almost surely if $j=0$ or
$i+j\geq4$.
The first case is trivial. We only need to consider $\bar S_{i,j}$ for
the cases $i=3, j=1$
and $i=2, j=2$, since the other ones are either contained or dual to these.
By Proposition~\ref{prop:skel}, $\bar S_{3,1}$ consists of
points which are almost surely in the skeleton and
where three paths coalesce. But the event that three coalescing
Brownian paths starting at distinct points coalesce at the same time
is almost surely empty.
By Proposition~\ref{prop:skel}, $\bar S_{2,2}$ consists of points
(almost surely in the double skeleton)
where two different forward paths coalesce and a backward path passes.
Since for any two forward and one backward Brownian
paths in the double skeleton, the event that this happens is almost
surely empty, by the perfectly coalescing/reflecting property of the
paths in the double skeleton
(see Subsubsection 3.1.1 and Theorem 8 of~\cite{kn:STW})
the conclusion follows.

Now, for the types that do occur.

{\bf Type (2,1) } By the above,
$S_{2,1}=\bar S_{2,1}$ almost surely, and  $\bar S_{2,1}$
consists almost surely of {\em points of coalescence}, that is
all points where two paths coalesce. By Proposition~\ref{prop:skel},
it is almost surely a subset of the skeleton, and thus is countable
(since there is at most one coalescence point for each pair of
paths starting from $\cal D$ in the skeleton). It is easy to see
that it is dense since the paths from a pair of nearby points in $\cal D$
also coalesce nearby with probability close to one.

{\bf Type (1,2) }
By the above, $S_{1,2}=\bar S_{1,2}$ almost surely, and
$\bar S_{1,2}$ consists almost surely of points where forward paths meet
backward paths. Thus, it is a subset of the (union of the traces of all the paths in the)
skeleton. It is easy to see that it is almost surely nonempty (and also dense).
We need only consider two such paths, say $W$ and $W^b$,
the former a forward one starting at $(0,0)$ (without loss of generality, by the
translation invariance of the law of $\W^D$), and the latter a backward one starting at an
arbitrary deterministic $(x_0,t_0)$, with $t_0>0$ to avoid a trivial case. It is clear that
the random set $\Lambda$ of space-time points $(t,W(t))$ for times $t\in[0,t_0]$ when
$W(t)=W^b(t)$ has a positive, less than
one probability of being empty. We will argue next the following claim.

{\noindent\bf Claim}
{\it $\Lambda$ has Hausdorff dimension $1$ for almost every
 pair of trajectories $(W,W^b)$
for which it is nonempty.}

By Proposition~\ref{prop:dpaths}, the distribution of $\{(W(t),W^b(t)):\,0\leq t\leq t_0\}$
(which is all that matters for this) can be
described in terms of two (forward) independent standard Brownian motions $B,B^b$ as follows
(see equations~(\ref{eq:for})-(\ref{eq:refl})). Let
$W^b(t)=x_0+B^b(t_0-t),\,t\leq t_0$, and  $\tau=\inf\{t\in[0,t_0]:\,B(t)=W^b(t)\}$,
with $\inf\emptyset=\infty$. If $\tau=\infty$, then $W=B$; otherwise, $W(t)=B(t)$
for $0\leq t\leq\tau$, and for $\tau\leq t\leq t_0$,
\beqnn
  W(t)=
  \begin{cases}
B(t)+\sup_{0\leq s\leq t}(W^b(s)-B(s)),
          \mbox{ if } W^b(0)<0;\\
B(t)-\sup_{0\leq s\leq t}(B(s)-W^b(s)),
          \mbox{ if } W^b(0)>0.
  \end{cases}
\eeqnn
Rewriting in terms of $W'(t):=W^b(t)-W^b(0),\,0\leq t\leq t_0$, which is a standard
Brownian motion independent of $B$, we have (for $0\leq t\leq t_0$)
$$
  W(t)=
  \begin{cases}
B(t)+\sup_{0\leq s\leq t}\{W'(s)-B(s)\}-W'(t_0)+x_0,
          \mbox{ if } W'(t_0)>x_0;\\
B(t)+\inf_{0\leq s\leq t}\{W'(s)-B(s)\}-W'(t_0)+x_0,
          \mbox{ if } W'(t_0)<x_0,
  \end{cases}
$$
if $\tau\leq t\leq t_0$, with $\tau=\inf\{t\in[0,t_0]:\,B(t)=W'(t)-W'(t_0)+x_0\}$;
otherwise, $W(t)=B(t)$.

 From the above discussion, we conclude that $\Lambda$
has the same distribution as the random set ${\cal G}$ obtained as follows. Let
${\cal T}^+$ and ${\cal T}^-$ be the sets of positive and negative record times of the
standard Brownian motion $X(t):=(W'(t)-B(t))/\sqrt2$, respectively,
i.e., ${\cal T}^+$ is the set of $t\geq0$ such that $X(t)=\sup_{0\leq s\leq t}X(s)$
and ${\cal T}^-$ is the same except with inf in place of sup. Consider also the
standard Brownian motion $Y(t):=(W'(t)+B(t))/\sqrt2$, which is independent of $X$.
If $W'(t_0)>x_0$, then
${\cal G}=\{([(X(t)+Y(t))/\sqrt2]-[(X(t_0)+Y(t_0))\sqrt2]+x_0,t):\,t\in{\cal T}^+\cap[\tau,t_0]\}$;
if $W'(t_0)<x_0$, then
${\cal G}=\{([(X(t)+Y(t))/\sqrt2]-[(X(t_0)+Y(t_0))\sqrt2]+x_0,t):\,t\in{\cal T}^-\cap[\tau,t_0]\}$.

It follows from Proposition~\ref{prop:haus} in Appendix~\ref{app:haus} that the sets
${\cal G}^\pm:=\{(X(t)+Y(t),t):\,t\in{\cal T}^\pm\cap[0,t_0]\}$ (one for each sign, respectively)
both have Hausdorff dimension $1$ almost surely.
Since the events $\{W'(t_0)>x_0\}$, $\{W'(t_0)<x_0\}$ and $\{\tau<t_0\}$ all have
positive probability, the claim follows.

{\bf Type (1,1) } $\bar S_{1,1}$ almost surely consists of
{\em points of continuation} of paths, that is,
all points
$(x,t)$ such that there is a path starting earlier than $t$ that touches
$(x,t)$.
By Proposition~\ref{prop:skel}, $\bar S_{1,1}$ is almost surely a subset
of the skeleton.
Since the trace of any single path has Hausdorff dimension $3/2$~\cite{kn:Ta}
and the countable union of such sets has the same dimension, it follows that
$\bar S_{1,1}$ has Hausdorff dimension $3/2$ almost surely.
By the previous parts of the proof, $\bar S_{1,1} \setminus S_{1,1}$
has lower dimension and so $S_{1,1}$ has the same Hausdorff dimension of $3/2$.

{\bf Type (0,1) } We claim that any
deterministic point is a.s.~of this type, hence
(by applying Fubini's Theorem)
$S_{0,1}$ is a.s.~of full Lebesgue measure in the plane.
That $m_{\mbox{{\scriptsize in}}}(x_0,t_0)=0$ a.s.~for every deterministic
$(x_0,t_0)$ follows from
Proposition~\ref{prop:skel}, since, if
$m_{\mbox{{\scriptsize in}}}(x_0,t_0)\geq1$, then there would be a
path in the skeleton passing through $(x_0,t_0)$, but this event clearly has
probability zero. The assertion that $m_{\mbox{{\scriptsize out}}}(x_0,t_0)=1$
a.s.~for every deterministic $(x_0,t_0)$ is
property $(o)$ of Theorem~\ref{teo:char}.

By Proposition~\ref{prop:dual}, the remaining types $(0,2)$ and $(0,3)$
are dual respectively to
$(1,1)$ and $(2,1)$, since the other types are dual to these.
Since ${\bar \W}^b$ is distributed like the standard Brownian web
(modulo a time reflection), the claimed results for types $(0,2)$ and $(0,3)$
follow from those already proved for $(1,1)$ and $(2,1)$.

\noindent{\bf Proof of Theorem~\ref{teo:typesb}}

{\bf Type (0,1) }
$\bar S_{1,1}$ is almost surely in the skeleton, thus making
$\bar S_{1,1}\cap\r\times\{t\}$ countable for all $t$. By a duality argument,
the same is true for $\bar S_{0,2}$. Since $\bar S_{0,1}=\r^2$ a.s.~by
Theorem~\ref{teo:types}, it follows that a.s.~for
all $t$, $S_{0,1}\cap\r\times\{t\}$ is of full Lebesgue
measure in the line.

Again, of the remaining types, it is enough by duality to
consider $(1,1)$, $(2,1)$ and $(1,2)$.

{\bf Type (2,1) }
For any deterministic $t$ and $(x_i,t_i)$ with $t_i<t$, $i=1,2$,
the probability that two independent Brownian paths starting
at $(x_i,t_i)$, $i=1,2$, respectively, coalesce exactly at time $t$
is zero. Since $S_{2,1}$ is in the skeleton,
$S_{2,1}\cap\r\times\{t\}=\emptyset$ almost surely.
Now, for any $t$, $|S_{2,1}\cap\r\times\{t\}|>1$ implies that
there are four independent Brownian paths starting at different
points, and such that the coalescence time of the first two and
that of the last two are the same. That this has zero probability
implies that a.s.~for all $t$, $|S_{2,1}\cap\r\times\{t\}|\leq 1$.

{\bf Type (1,2) }
For any deterministic $t$,
$S_{1,2}\cap\r\times\{t\}=\emptyset$ almost surely, since the
probability of two fixed paths, one forward, one backward, meeting
at a given deterministic time is $0$.
Indeed, from the analysis of type $(1,2)$ done above in the proof of
Theorem~\ref{teo:typesa}, this is because the probability that a
Brownian motion has a record value at a given deterministic time
is~$0$.
For any $t$, $|S_{1,2}\cap\r\times\{t\}|>1$ implies that
there exist in the double Brownian web skeleton two pairs,
each consisting of
one forward and one backward path, such that in both pairs
the forward and backward paths meet at the same time.
We claim that this has zero probability
and thus that $|S_{1,2}\cap\r\times\{t\}|\leq 1$ almost surely.
To verify the claim, we again use the analysis of type $(1,2)$ done
for Theorem~\ref{teo:typesa}, which shows that it suffices to
prove that there is zero probability that
two independent standard Brownian motions $B_1,\,B_2$
have a common strictly positive record time. But, as noted
in Appendix~\ref{app:haus}, this is the same as having zero
probability for $B_1,\,B_2$ to both have a zero at
a common strictly positive time. This latter
probability is indeed zero because of
the well known fact that the two-dimensional Brownian motion
$(B_1,B_2)$ a.s.~does not return to $(0,0)$.

{\bf Type (1,1) } Since points with
$m_{\mbox{{\scriptsize in}}}\geq1$ are a.s.~in the skeleton,
$\bar S_{1,1}\cap\r\times\{t\}$ is a.s.~countable
(and easily seen to be dense) for every
$t\in\r$. Now the previous parts of the proof
imply that the same holds for
$S_{1,1}\cap\r\times\{t\}$ for every $t\in\r$.

%%%%%%%%%%%%%%%%%%%%%%%%%%%%%%%%%%%%%%%%%%%%%%%
%%%%%%%%%%%%% SECTION 6 %%%%%%%%%%%%%%%%%%%%%%%%%%%%
%%%%%%%%%%%%%%%%%%%%%%%%%%%%%%%%%%%%%%%%%%%%%%%
%%%%%%%%%%%%%%%%%%%%%%%%%%%%%%%%%%%%%%%%%%%%%%%

\section{General convergence results}
\setcounter{equation}{0}
\label{sec:conv}

In this section, we state and prove Theorem~\ref{thm:gen},
which is an extension
of our convergence result
for noncrossing paths, Theorem~\ref{n-c-conv},
to the case where paths can cross (before the scaling
limit has been taken). At the end of the section, we show
that the noncrossing Theorem~\ref{n-c-conv} follows
from Theorem~\ref{thm:gen} and other results.

Before stating our general theorem that allows
crossing, we briefly
discuss some systems with crossing paths, to which it
should be applicable. Consider the stochastic difference
equation~(\ref{diffeq}) where the $\Delta_{i,j}$'s are
i.i.d.~integer-valued random variables, with zero mean
and finite nonzero variance. Allowing $(i,j)$ to be
arbitrary in $\Z ^2$, we obtain as a natural generalization
of Figure 1, a collection of random piecewise linear paths
that can cross each other, but that still coalesce when they meet at
a lattice point in $\Z ^2$.

With the natural choice of diffusive
space-time scaling and under
an irreducibility condition (to insure that the walks
from any two starting points have a strictly positive probability
of coalescing), the scaling limit of such a discrete
time system should be the standard Brownian web. To
see what happens in reducible cases, consider simple
random walks ($\Delta_{i,j} = \pm 1$), where the paths
on the even and odd subsets of $\Z ^2$ are independent
of each other, and so the scaling limit on all of $\Z ^2$
consists of the union of two independent Brownian webs.
For $\Delta_{i,j} = \pm 2$, the limit would be the union
of four independent Brownian webs. We remark that for
continuous time random walks (as discussed in the next
section of this paper for $\Delta_{i,j} = \pm 1$),
no irreducibility condition is
needed.

We proceed with some definitions needed for our general
convergence theorem.
For $a,b,t_0\in\r,\, a<b$, and $t>0$, we define two real-valued measurable
functions $l_{t_0,t}([a,b])$
and $r_{t_0,t}([a,b])$ on $(\h,\f_\h)$ as follows.
For $K \in \h$ ,
$l_{t_0,t}([a,b])$ evaluated at $K$ is defined as
  $ \inf\{x \in [a,b]|  \exists y \in \R$
and  a path in $ K$  which touches both $ (x,t_0) \hbox{ and }
(y,t_0 + t)\}$ and
$r_{t_0,t}([a,b])$ is defined similarly with the $\inf$ replaced by
   $ \sup$. We also define the following functions on $(\h,\f_\h)$
whose values are subsets of $\R$. As before we let $K \in \h$ and
suppress
$K$ on the left hand side of the formula for ease of notation.
\begin{eqnarray}
\nonumber
{\cal N}_{t_0,t} ([a,b])
& =& \{ y \in \R\,|\,\exists\, x \in [a,b]
\hbox{ and a path in } K \hbox{ which touches} \\
  \label{eq:caln}
&&\hspace{1.6cm}\hbox{both }  (x,t_0) \hbox{ and }
(y,t_0 + t)\}
\end{eqnarray}
\begin{eqnarray}
\nonumber
{\cal N}^-_{t_0,t} ([a,b])
& =& \{ y \in \R\,|\, \hbox{
  there is a path in } K
\hbox{ which touches both }\\
  \label{eq:caln-}
&&\hspace{1.6cm} (l_{t_0,t}([a,b]),t_0) \hbox{ and }
(y,t_0 + t)\}
\end{eqnarray}
\begin{eqnarray}
\nonumber
{\cal N}^+_{t_0,t} ([a,b])
& =& \{ y \in \R\,|\, \hbox{
  there is a path in } K\hbox{ which touches both } \\
  \label{eq:caln+}
&&\hspace{1.6cm} (r_{t_0,t}([a,b]),t_0) \hbox{ and }
(y,t_0 + t)\}
\end{eqnarray}

\brm
We notice that $|{\cal N}_{t_0,t} ([a,b])|= \eta(t_0,t;a,b)$.
\erm

Let $\{\X_m\}$ be a sequence of $(\h,\f_\h)$-valued random variables
with distributions $\{\mu_m\}$. We define conditions $(B'_1), (B'_2)$
as follows.

\noindent
\beqnn
(B'_1)&&\!\!\!\!\!\!\!\!  \forall\beta >0,\,\,
\limsup_{m\to \infty}\sup_{t >\beta}  \sup_{t_0,a \in \R} \mu_m (
|{\cal N}_{t_0,t} ([a-\e,a+\epsilon]) | > 1) \to 0 \hbox{ as } \epsilon
\to 0^+\,\\\nn
(B'_2)&&\!\!\!\!\!\!\!\!  \forall\beta >0,
\frac{1}{\e} \limsup_{m\to \infty}\sup_{t >\beta}
\sup_{t_0,a \in \R} \mu_m (
{\cal N}_{t_0,t} ([a-\e,a+\epsilon])\ne {\cal N}^+_{t_0,t}([a-\e,a+\e])\\
&&\hspace{3.8cm}
\cup{\cal N}^-_{t_0,t}([a-\e,a+\e]))\to0\hbox{ as }\epsilon\to0^+
\eeqnn

\brm
\label{mono}
Note that if we consider a process with non-crossing paths then
conditions $(B'_1)$ and $(B'_2)$ follow from conditions $(B_1)$ and $(B_2)$
respectively because of the following monotonicity property.
For all $a < b, t_0$ and $0 < s < t$
\[
\P(|\eta(t_0,t;a,b)| \geq k) \leq \P(|\eta(t_0,s;a,b)| \geq k)
\]
for all $k \in \N$.
\erm

\bteo
\label{thm:gen}
Suppose that $\{ \mu_m\}$ is tight. If
Conditions $(I_1), (B'_1) $ and $(B'_2)$ hold,
then $\{\X_m\}$ converges in distribution to the Brownian web $\bw$.
\eteo

\brm
\label{rmk:gcdbw}
There is a natural analogue to (and corollary of) this theorem for the
double Brownian web, in which ($I_1$) is replaced by Condition ($I_1^D$)
of Theorem~\ref{teo:convdb} below and ($B'_1$), ($B'_2$) and their
backward analogues are separately valid. The proof of this corollary
is like that of Theorem~\ref{teo:convdb}, except based on the general
convergence theorem for the (forward) Brownian web rather than the special
convergence theorem for the noncrossing case.
\erm

Theorem~\ref{thm:gen} is proved through a series of lemmas.

\blem
\label{lem:A}
Let $\mu$ be a subsequential limit of $\{ \mu_m\}$ and suppose that
$\mu$ satisfies condition $(i')$ of Theorem $~\ref{teo:char2}$ and
\beqnn
(B''_1)&& \!\!\!  \forall\beta >0,\,\,\,\sup_{t >\beta} \sup_{t_0,a} \mu(| {\cal N}_{t_0,t}
([a-\e,a+\epsilon])| > 1)
\to 0 \hbox{ as } \e \to 0^+\\\nn
(B''_2)&&\!\!\!  \forall\beta >0,\,\frac{1}{\e} \sup_{t >\beta} \sup_{t_0,a}
\mu({\cal N}_{t_0,t} ([a-\e,a+\epsilon])
\ne {\cal N}^+_{t_0,t}
([a-\e,a+\e])\\
&&\hspace{2.2cm}
\cup{\cal N}^-_{t_0,t} ([a-\e,a+\e]))\to 0\hbox{ as }\epsilon\to0^+.
\eeqnn
\noindent Then $\mu$ is the distribution of the Brownian web.
\elem
{\bf Proof }
It follows from condition $(i')$ and $(B''_1)$ that the limiting
random variable $\X$ satisfies
condition $(i)$ of the characterization theorem~\ref{teo:char}. That is ($\mu$)
almost surely there is exactly one path starting from
each point of $\cd$ and these paths are distributed as coalescing
Brownian motions. Let us define a $(\h,\f_\h)$-valued random variable
$\X'$ on the same probability space as the one on which $\X$ is defined
to be
the closure in $(\Pi,d)$ of the paths of $\X$ starting from $\cd$. We will
denote
probabilities in the common probability space by $\P$.
$\X'$ has the distribution of $\bw$.
We need to show that it also satisfies
condition $(ii''')$
of Theorem $~\ref{teo:char2}$.
\noindent Let $a < b, t_0 \in \R $ and $t > 0$ be given. For the
random variable $\X$ we will denote the
counting random variable $\eta
(t_0,t
;a,b)$ by $\eta$ and the corresponding variable for $\X'$ by $\eta'$. Let
$z_j = (a + j (b - a)/M,t_0)$ for $j = 0,1, \ldots,M$, be $M+1$ equally
spaced points in the interval $[a,b]$.

Now define
\noindent $
\eta_M = |\{x \in \R| \exists \hbox{ a path in $\X$
which
touches both a point in }$

\noindent $\{z_0,\ldots,z_M \} \hbox{ and } (x,t+t_0)\}|
$, where $|\cdot|$ stands for cardinality.
 Let $\eta'_M$
be the corresponding random variable for $\X'$.
Clearly $\eta \geq \eta_M \hbox{ and } \eta' \geq \eta'_M $.
 From ($B_1''$) it follows that $\eta_M = \eta'_M$ almost
surely.
Now let $\e = \frac{(b-a)}{M}$. By condition $(B''_2)$,
letting $M \to \infty\, (\e \to 0)$, we obtain
\[
\P(\eta > \eta'_M) = \P(\eta > \eta_M) \to 0 \hbox{ as } M \to \infty.
\]
Thus, $\P(\eta > \eta') = 0$, showing that $\eta$ is stochastically
dominated
by $\eta'$. This completes the proof of the lemma.

For $t > 0,\, \e>0,\,0< \e' < \frac{\e}{8}, \,0 \leq \d < \frac{t}2$,
consider  the following event.
\begin{eqnarray*}
&O(a,t_0,t,\e, \e', \d)=&\\
&\{K \in \h| \mbox{ there are three paths }
(x_1(t),t_1),(x_2(t),t_2),(x_3(t),t_3)
\mbox{ in }K&\\
&\mbox{with } t_1,t_2,t_3 <t_0+ \d,
x_1(t_0+\d) \in (a - \e-\e',a-\e+\e'),&\\
&x_2(t_0+\d) \in (a-\e+2\e',a+\e-2\e'),
x_3(t_0 + \d) \in (a+\e-\e',a+\e+\e')&\\
&\mbox{and } x_2(t_0+t) \ne x_1(t_0 +t),
x_2(t_0+t) \ne x_3(t_0 +t)\}&
 \end{eqnarray*}

\blem
\label{lem:B}
$(B''_2)$ in Lemma ~$\ref{lem:A}$ can be replaced by:
\[
(B_2''')\,\,\, \forall\beta>0,\,\frac{1}{\e} \limsup_{\e' \to 0}\sup_{t >\beta}\sup_{t_0,a}
\limsup_{\d \to 0}
  \mu (O(a,t_0,t,\e, \e', \d)) \to 0
\hbox{ as } \e \to 0^+.
\]
\elem

{\bf Proof } We prove the lemma by showing that conditions $(i')$ and
$(B''_1)$ together with $(B_2''')$ imply condition
$B''_2$. Let $\beta > 0$. Define $C_1(b,t_0,\e',\d)$ as
 \begin{eqnarray*}
\{\hspace{-.7cm}&&
K \in \h|\hbox{ there is a path in } K \hbox{ which touches both }
(b,t_0)\\
&&\hspace{1.4cm} \hbox{ and } \{b -\e'\} \times [t_0,t_0+\d]
\cup\{b + \e'\} \times [t_0,t_0+\d] \},
 \end{eqnarray*}
\noindent and
$C_2(a,t_0,\e,\e',\d)$ as
 \begin{eqnarray*}
\{\hspace{-.7cm}&&
K \in \h|\hbox{ there is a path in } K \hbox{ which touches both }
[a-\e,a+\e] \times \{t_0\}\\
&&\hspace{1.4cm}  \hbox{ and } \{a -\e-\e'\} \times
[t_0,t_0+\d]\cup
\{a + \e+\e'\} \times [t_0,t_0+\d] \}.
 \end{eqnarray*}
 \noindent Now observe that (modulo sets of zero $\mu$ measure)
 \begin{eqnarray*}
&\{{\cal N}_{t_0,t} ([a-\e,a+\e])
\ne{\cal N}^+_{t_0,t}([a-\e,a+\e])\cup{\cal N}^-_{t_0,t}([a-\e,a+\e])\}&\\
&\cap C^c_1(a+\e,t_0,\e',\d)\cap C^c_1(a-\e,t_0,\e',\d)
\cap C^c_2(a,t_0,\e,\e',\d)&\\
&\cap
\{|{\cal N}_{t_0+\d,t-\d}([a-\e-2\e',a-\e+2\e'])| =1\}&\\
&\cap
\{|{\cal N}_{t_0+\d,t-\d} ([a+\e-2\e',a+\e+2\e'])| =1\}&\\
&\subseteq O(a,t_0,t,\e,\e',\d).&
\end{eqnarray*}
\noindent Therefore, we have
\begin{eqnarray*}
&&\!\!\!\mu({\cal N}_{t_0,t} ([a-\e,a+\e])  \ne {\cal N}^+_{t_0,t}
([a-\e,a+\e])
\cup  {\cal N}^-_{t_0,t} ([a-\e,a+\e])) \\
&\leq&\!\!\!\mu(O(a,t_0,t,\e,\e',\d))
+\mu(C_2(a,t_0,\e,\e',\d))
  + \mu(C_1(a+\e,t_0,\e',\d)) \\
&+&\!\!\!\mu(C_1(a-\e,t_0,\e',\d))
  +\mu( |{\cal N}_{t_0+\d,t-\d} ([a-\e-2\e',a-\e+2\e'])| >1)\\
&+&\!\!\!\mu(|{\cal N}_{t_0+\d,t-\d} ([a+\e-2\e',a+\e+2\e'])| >1).
\end{eqnarray*}
Letting $\d \to 0$, we obtain
\begin{eqnarray}
&&\!\!\!\mu ({\cal N}_{t_0,t} ([a-\e,a+\e])
\ne {\cal N}^+_{t_0,t}
([a-\e,a+\e])
\cup  {\cal N}^-_{t_0,t} ([a-\e,a+\e])) \leq \nonumber \\
&&\!\!\!\limsup_{\d \to 0}\{\mu(O(a,t_0,t,\e,\e',\d))
+\mu(C_2(a,t_0,\e,\e',\d))
  + \mu(C_1(a+\e,t_0,\e',\d)) \nonumber \\
&+&\!\!\!  \mu(C_1(a-\e,t_0,\e',\d))
+\mu( |{\cal N}_{t_0+\d,t-\d} ([a-\e-2\e',a-\e+2\e'])| >1) \nonumber \\
&+&\!\!\!\mu(|{\cal N}_{t_0+\d,t-\d} ([a+\e-2\e',a+\e+2\e'])| >1)\}.\label{open}
\end{eqnarray}
Now,
\beqnn
\lim_{\d \to 0}  \mu(C_1(a+\e,t_0,\e',\d)) &=&
\lim_{\d \to 0} \mu(C_1(a-\e,t_0,\e',\d))\\ &=& \lim_{\d \to 0}
\mu(C_2(a,t_0,\e,\e',\d)) = 0,
\eeqnn
\noindent since elements of $\h$ are compact subsets of $\Pi$,
and paths in $\Pi$ cannot have
completely flat segments.

Now since, $t - \d >\frac{t}{2}> \frac{\beta}{2}$, it follows from
$(B''_1)$ that
for all $\gamma >0$,
\beqnn
& &\sup_{t>\beta}\sup_{a,t_0}\sup_{0<\d<\frac{t}{2}
 }\mu(|{\cal N}_{t_0+\d,t-\d}
([a-\gamma,a+ \gamma])| >1| ) \\
& \leq &
\sup_{t > \frac{\beta}{2}}\sup_{a,t_0}\mu(|{\cal N}_{t_0,t}
([a-\gamma,a+ \gamma])| >1| ) \to 0 \hbox{ as } \gamma \to 0. \\
\eeqnn
This implies that for all $\e >0$,
\beq
\limsup_{\e' \to 0} \sup_{t>\beta}\sup_{a,t_0}
\limsup_{\d \to 0}
\mu( |{\cal N}_{t_0+\d,t-\d} ([a\pm\e-2\e',a\pm\e+2\e'])| >1)
 =0.
\eeq
Together with~(\ref{open}), this gives us
\begin{eqnarray}\nn
\sup_{t>\beta}\sup_{t_0,a}\mu ({\cal N}_{t_0,t}
([a-\e,a+\e])  \ne {\cal N}^+_{t_0,t}
([a-\e,a+\e])
\cup  {\cal N}^-_{t_0,t} ([a-\e,a+\e]))\\
\leq \limsup_{\e' \to 0} \sup_{t>\beta} \sup_{a,t_0,} \limsup
_{\d \to 0}
\mu(O(a,t_0,t,\e,\e',\d))\hspace{3cm}
\end{eqnarray}
Now, using $(B_2''')$, we obtain
\[
\frac{1}{\e}
\sup_{t>\beta}\sup_{t_0,a}
\mu ({\cal N}_{t_0,t}([a-\e,a+\e])\ne{\cal N}^+_{t_0,t}([a-\e,a+\e])
\cup{\cal N}^-_{t_0,t} ([a-\e,a+\e]))\to0
\]
as $\e \to 0^+$, proving the lemma.

{\bf Proof of Theorem~\ref{thm:gen} } Tightness implies that every
subsequence of $\{\mu_m\}$ has
a subsequence converging to some $\mu$. Let us denote
the corresponding limiting random variable by $\X$. We prove the theorem by
showing that every such
$\mu = \mu_{\bw}$. From Lemmas $\ref{lem:A}$ and $\ref{lem:B}$ it
follows that it is sufficient to prove condition $(i')$ of
Theorem~\ref{teo:char2}, condition $(B''_1)$ and condition $(B'''_2)$.

Let  $\beta>0$ and define for all  $0 \leq \d < t/2$,
$
{\cal N'}^\d_{t_0,t} ([a,b]) = \{ y \in \R |\exists \hbox{ a path }$

\noindent $(x(s),s_0), s_0 < t_0 + \d
\hbox{ in } K  \hbox{ such that }  x(t_0+\d) \in (a,b) \hbox{ and }
x(t_0 + t)= y \}$.
We note that
the set $\{|{\cal N'}^\d_{t_0,t} ([a,b])| > 1 \}$ is an open subset of
$\h$
for all $\d \geq 0$.
Then we have
\beqnn
& &\sup_{t>\beta}\sup_{t_0,a}\mu(|{\cal N}_{t_0 ,t}
([a-\e,a + \e])| > 1 )\\ 
 &\leq&
\sup_{t>\beta}\sup_{t_0,a}
\limsup_{\d \to 0}\{\mu(|{\cal N'}^\d_{t_0,t}
([a-2 \e,a+2 \e])| > 1 )\\
&+&\mu(C_2(a,t_0,\e,\e,\d))\}\\
& \leq &\sup_{t>\frac{\beta}{2}}\sup_{t_0,a}\mu(|{\cal N'}^0_{t_0,t}
([a-2 \e,a+2 \e])| > 1 )\\
&+ & \sup_{t_0,a} \limsup_{\d \to 0} \mu(C_2(a,t_0,\e,\e,\d))
\eeqnn
Now, $\limsup_{\d \to 0}\mu(C_2(a,t_0,\e,\e,\d)) = 0 $, since elements
of $\h$ are compact subsets of $\Pi$. This together with
the fact that $\{|{\cal N'}^\d_{t_0,t} ([a,b])| > 1 \}$ is an open subset of
$\h$ leads to

\beqnn
&  &\sup_{t>\beta}\sup_{t_0,a} \mu(|{\cal N}_{t_0 ,t}
([a-\e,a + \e])| > 1 )  \\ 
  &\leq & \sup_{t>\frac{\beta}{2}}\sup_{t_0,a} \mu(|{\cal N'}^0_{t_0 ,t}
([a-2\e,a + 2\e])| > 1 )  \\
& \leq& \sup_{t>\frac{\beta}{2}}\sup_{t_0,a} \limsup_m
\mu_m(|{\cal N'}^0_{t_0 ,t}
([a-2\e,a + 2\e])| > 1 )\\
& \leq& \sup_{t>\frac{\beta}{2}}\sup_{t_0,a} \limsup_m
\mu_m(|{\cal N}_{t_0 ,t}
([a-2\e,a + 2\e])| > 1 )\\
& \leq &   \limsup_m \sup_{t>\frac{\beta}{2}}\sup_{t_0,a}
\mu_m(|{\cal N}_{t_0 ,t}
([a-2\e,a + 2\e])| > 1 ).
\eeqnn
It follows from $(B'_1)$ that
\[
\limsup_m \sup_{t>\frac{\beta}{2}}\sup_{t_0,a}
\mu_m(|{\cal N}_{t_0 ,t}
([a-2\e,a + 2\e])| > 1 ) \to 0 \hbox{ as } \e \to 0^+.
\]
This proves $(B''_1)$, which implies that
\begin{itemize}
\item[$(o)$] starting from any deterministic point, there is
$\mu$-almost  surely
only a single path in $\X$.
\end{itemize}
Combining this with $(I_1)$, we readily obtain that
\begin{itemize}
\item[$(i)$] the finite-dimensional distributions of $\X$ are those of
coalescing
Brownian motions with unit diffusion constant.
\end{itemize}
Condition
$(i')$ of Theorem
$~\ref{teo:char2}$ follows immediately from $(o)$ and $(i)$. Now we
proceed to verify condition $(B'''_2)$.

We have
\beqnn
&&\sup_{t>\beta}\sup_{t_0,a}
\limsup_{\d \to 0}
  \mu (O(a,t_0,t,\e, \e', \d)) \\
&\leq  &
\sup_{t>\frac{\beta}{2}}\sup_{a,t_0}
 \mu (O(a,t_0,t,\e, \e', 0))\\
&\leq &
\limsup_m \sup_{t>\frac{\beta}{2}}\sup_{a,t_0}
 \mu_m (O(a,t_0,t,\e, \e', 0))\\
&\leq &
\limsup_m \sup_{t>\frac{\beta}{2}}\sup_{a,t_0}
\mu_m ({\cal N}_{t_0,t}
([a-\e-\e',a+\e+\e']) \\
& & \ne {\cal N}^+_{t_0,t}
([a-\e-\e',a+\e+\e'])
\cup  {\cal N}^-_{t_0,t} ([a-\e-\e',a+\e+\e'])),
\eeqnn
where the second inequality follows from the fact that
$O(a,t_0,t,\e, \e', 0)$ is an open subset of $\h$.
For the third inequality to hold we need to insure that there is
no more than
one path touching either $(a-\e-\e',t_0) \hbox{ or } (a+\e+\e',t_0)$;
this follows from $(B'_1)$.
Since $\e' < \frac{\e}{8},\,\e + \e' \to 0 \hbox{ as } \e \to 0 $,
using condition $(B'_2)$, we obtain
$$\frac{1}{\e} \limsup_{{\e'} \to 0}\sup_{t_0,a}
\limsup_{\d \to 0}
  \mu (O(a,t_0,t,\e, \e' ,\d)) \to 0 \hbox{ as } \e \to 0^+,$$
proving condition $(B'''_2)$. This completes the proof
of the theorem.

We now suppose that $\X_1, \X_2,\ldots$ is a sequence of
$({\cal H},{\cal F}_{{\cal H}})$-valued random variables
so that each $\X_i$ consists of {\em noncrossing} paths.
The noncrossing condition produces a considerable simplification
of Theorem 5.1; namely Theorem~\ref{n-c-conv}.

\noindent{\bf  Proof of Theorem \ref{n-c-conv}} This is an immediate
consequence of Remark~\ref{mono},
 Theorem \ref{thm:gen} and Proposition \ref{prop: nctightness}
of the appendix.

%%%%%%%%%%%%%%%%%%%%%%%%%%%%%%%%%%%%%%%%%%%%%%
%%%%%%%%%%%%% SECTION 7 %%%%%%%%%%%%%%%%%%%%%%%%%%%
%%%%%%%%%%%%%%%%%%%%%%%%%%%%%%%%%%%%%%%%%%%%%%
%%%%%%%%%%%%%%%%%%%%%%%%%%%%%%%%%%%%%%%%%%%%%%

\section{Convergence for coalescing random walks}

\setcounter{equation}{0}
\label{sec:rw-conv}

We now apply Theorem~\ref{n-c-conv} to coalescing random walks.
For that, we begin by precisely defining $Y$ (resp., $\tilde Y$),
the set of all discrete (resp., continuous) time coalescing random walks
on $\Z$. For $\delta$ an arbitrary positive real number, we obtain sets
of rescaled walks,
$Y^{(\d)}$ and $\tilde Y^{(\d)}$, by the usual rescaling of space by
$\d$ and time by $\d^2$.
The (main) paths of $Y$ are the discrete-time random walks $\yo$, as
described in the Introduction
and shown in Figure~1, with $(y_0,s_0)= (i_0,j_0) \in\Z \times \Z$
arbitrary except that $i_0+j_0$
must be even. Each random walk path goes from $( i, j)$ to $(i\pm1,j+1)$
linearly. In addition to these, we add some boundary paths so that $Y$
will be a compact subset of $\o$. These are all the paths of the form
$(f,s_0)$ with $s_0\in \Z \cup \{-\infty,\infty\}$ and $f\equiv\infty$
or
$f\equiv-\infty$. Note that for $s_0=-\infty$ there are two different
paths starting from the single point at $s_0=-\infty$ in $\br^2$.

The continuous time $\tilde Y$ can be defined similarly,
except that here $y_0$ is any $i_0\in\Z$ and $s_0$ is arbitrary
in $\R$. Continuous time walks are normally seen as jumping from
$ i$ to $i\pm1$
at the times $T^{(i)}_k\in(-\infty,\infty)$ of a rate one Poisson
process. If the jump is, say, to $i+1$, then our polygonal path will
have a linear segment between $(i,T^{(i)}_k)$ and
$(i+1,T^{(i+1)}_{k'})$, where $T^{(i+1)}_{k'}$ is the first Poisson
event at $i+1$ after $T^{(i)}_k$.
Furthermore, if $T^{(i_0)}_{k}<s_0<T^{(i_0)}_{k+1}$, then there will be
a constant segment in the path before the first nonconstant linear
segment. If $s_0=T^{(i_0)}_{k}$, then we take two paths: one with an
initial constant segment and one without.

\bteo
\label{teo:convrw}
Each of the collections of rescaled coalescing random
walk paths, $Y^{(\d)}$ (in discrete time) and
$\tilde Y^{(\d)}$ (in continuous time) converges in distribution
to the standard Brownian web as $\d \to 0$.
\eteo

\noindent{\bf Proof }
By Theorem \ref{n-c-conv}, it suffices to verify conditions
$I_1, B_1$ and $B_2$.

Condition $I_1$
is basically a consequence of the Donsker
invariance principle, as already noted in the Introduction.
Conditions $B_1$ and $B_2$ follow from the coalescing walks
version of the inequality of (\ref{etabound}), which is
\begin{equation}
\label{walkbound}
\mu_\d(\eta(t_0,t;a,a+\e)\geq k)\, \leq \,
[\mu_\d(\eta(t_0,t;a,a+\e)\geq2)]^{k-1}.
\end{equation}
Taking the sup over $(a,t_0)$ and the $\limsup$ over $\d$ and
using standard random walk arguments produces
an upper bound of the from $C_k (\e/ \sqrt{t})^{k-1}$
which yields $B_1$ and $B_2$ as desired.

We next state extensions of Theorems~\ref{n-c-conv}
and~\ref{teo:convrw} to the
Double Brownian web.

\bteo
\label{teo:convdb}
Suppose $\X_1^D=(\X_1,\X_1^b),\X_2^D=(\X_2,\X_2^b),\ldots$ are
$({\cal H}^D,{\cal F}_{{\cal H}^D})$-valued random variables
such that each $\X_n$ consists only of noncrossing paths and similarly
for each $\X_n^b$. Suppose further that there is a dense countable
deterministic
subset $\cal D$ of $\r^2$ such that
\begin{itemize}
\item[$(I_1^D)$] there exists $\theta^y_n\in \X_n$ and
${\theta_n^y}^b\in \X_n^b$
such that for any deterministic
$y_1,\ldots,y_m\in{\cal D},\{\theta_n^{y_1},{\theta_n^{y_1}}^b,
\ldots,\theta^{y_m}_n,{\theta_n^{y_m}}^b\}$ converge in distribution
to coalescing/reflecting forward/backward Brownian motions.
\end{itemize}
If in addition conditions ($B_1$) and  ($B_2$) of Theorem~\ref{n-c-conv}
above are valid for $\hat\eta$ and $\hat\eta^b$ (separately), then
$\X_n^D$
converges in distribution to the double Brownian
web.
\eteo

\noindent{\bf Proof } This could be proved as a corollary of a general
(i.e., allowing paths that cross) convergence to the double Brownian web
result (see Remark~\ref{rmk:gcdbw}); instead we prove it as a corollary to
Theorem~\ref{n-c-conv} as follows.
By Prop.~\ref{prop: nctightness},
($I_1^D$) implies
tightness separately for the forward and backward
processes; this, in turn, implies tightness of the joint process. This and
$(I_1^D)$ imply that any subsequential limit $\X^D=(\X,\X^b)$ of $(\X^D_n)$ contains
a version of the DBW, say $\Y^D=(\Y,\Y^b)$. But, from Theorem~\ref{n-c-conv},
separate convergence of $\X_n$ and
$\X^b_n$ to the forward and backward BW's, respectively,
implies that $\X$ and $\X^b$ are versions of the forward and backward BW's,
respectively. Now, Theorem~\ref{teo:charm} implies that
$\X=\Y$ and $\X^b=\Y^b$ (almost surely), and the argument is complete.

We now consider the joint convergence of forward coalescing random walks
and
the dual backward coalescing random walks.
The backward (dual) $Y^{b,(\delta)}$ and $\tilde Y^{b,(\delta)}$
processes for discrete and continuous time random walks can be defined in
a straightforward  way. In discrete time, they are polygonal paths on the dual
lattice of all $(i,j)$ with $i+j$ {\em odd}; in continuous time, they are
the paths of walkers on the dual lattice $\z+\frac12$ (see~\cite{kn:STW}
and~\cite{kn:FINS1,kn:FINS2} for more details).

\bteo
\label{teo:convdbl}
$({Y}^{(\delta)},\tilde Y^{b,(\delta)})$ and $(\tilde {Y}^{(\delta)},\tilde
Y^{b,(\delta)})$
converge in distribution as $\delta\to0$ to the double Brownian web.
\eteo

\noindent{\bf Proof }
The validity of $(I_1^D)$ is shown in Subsubsection 3.3.6
of~\cite{kn:STW}; $(B_1)$ and
$(B_2)$ for the forward and backward processes
separately have already been shown.

\vspace{.5cm}

\noindent{\bf Acknowledgements} The authors thank Raghu Varadhan for many 
useful discussions, Balint T\'oth for comments related to reference~\cite{kn:TW}, 
Yuval Peres for correspondence related to Appendix~\ref{app:haus}, 
and Anish Sarkar for discussions and references on drainage networks.
They also thank the Abdus Salam International Centre for Theoretical
Physics for its hospitality when an early draft of this work was
being prepared.
L.R.G.F.~thanks the Courant Institute of NYU and the Math
Department at the University of Rome "La Sapienza" for hospitality
and support during visits where parts of this work were done.
M.I.~thanks the Courant Institute of NYU and the IME at the University 
of S\~ao Paulo for hospitality and support during visits where parts of 
this work were done. K.R.~thanks the Courant Institute of NYU, 
where parts of this work were done,
for hospitality and support during a sabbatical leave visit.

\appendix

%%%%%%%%%%%%%%%%%%%%%%%%%%%%%%%%%%%%%%%%%%%%%%%
%%%%%%%%%%%%%%%% APPENDIX A %%%%%%%%%%%%%%%%%%%%%%%%
%%%%%%%%%%%%%%%%%%%%%%%%%%%%%%%%%%%%%%%%%%%%%%%
%%%%%%%%%%%%%%%%%%%%%%%%%%%%%%%%%%%%%%%%%%%%%%%

\section{Some measurability issues}
\label{app:mea}
\setcounter{equation}{0}

Let $({\cal H}, d_{\cal H})$ denote the Hausdorff metric space
induced by $(\Pi,d)$. $\f_\h$ denotes the
$\sigma$--field
generated by the open sets of  $\cal H$. We will consider now
{\em cylinders} of $\cal H$.
Let us fix nonempty horizontal segments
$I_1,\ldots,I_n$ in $\R^2$ (i.e., $I_k=I'_k\times\{t_k\}$),
where each $I'_k$ is an interval (which need not be finite and can be
open, closed or neither)
and $t_k\in\R$. Define
\beqn
\nonumber
&&C_{I_1,\ldots,I_n}^{t_0}:=\{K\in{\cal H}:
\mbox{ there exists $(f,t)\in K$ such that}\\
\label{eq:cil}
&&\quad\quad\quad\quad\quad\,\,
t>t_0 \ \mbox{and $(f,t)$ goes through $I_1,\ldots,I_n$}\}\\
\nonumber
&&\bar C_{I_1,\ldots,I_n}^{t_0}:=\{K\in{\cal H}:
\mbox{ there exists $(f,t)\in K$ such that}\\
\label{eq:cil1}
&&\quad\quad\quad\quad\quad\,\,
t\geq t_0 \ \mbox{and $(f,t)$ goes through $I_1,\ldots,I_n$}\}\\
\nonumber
&&C_{I_1,\ldots,I_n}:=\{K\in{\cal H}:
\mbox{ there exists $(f,t)\in K$ such that}\\
\label{eq:cil2}
&&\quad\quad\quad\quad\quad\,\,
(f,t)\mbox{ goes through $I_1,\ldots,I_n$}\}.
\eeqn

We will call sets of the form~(\ref{eq:cil}) {\em open cylinders} if
each $I_k$ is open, and sets of the of the form~(\ref{eq:cil1}) {\em closed cylinders} if
each $I_k$ is closed.

\brm
\label{sigid}
It is easy to see that sets of the
form~(\ref{eq:cil}-\ref{eq:cil2}) for arbitrary $I_1,\ldots,I_n$
can be generated by open cylinders.
\erm

Let now
$\mathfrak C$ be the $\sigma$--field generated by the open cylinders.
\bprop
\label{sigid1}
$\f_\h=\mathfrak C$
\eprop

The proposition is a consequence of the following two lemmas.

\blem
\label{sigid2}
$\f_\h\supset\mathfrak C$
\elem

\noindent{\bf Proof }
It is enough to observe that the open cylinders are open sets
of $\cal H$. Indeed, take an open cylinder, an element $K$ in that cylinder,
and $(f,t)\in K$ such that $t>t_0$ and
$a_i<f(t_i)<b_i$ for all $i=1,\ldots,n$.
All points of $B_{\cal H}(K,\epsilon)$,
the open ball in $\cal H$ around $K$ with radius $\epsilon$,
contains a path $(f',t')$ in a ball in $\Pi$ around $(f,t)$ of
radius $\epsilon$. Thus by choosing $\epsilon$ small enough,
$(f',t')$ will satisfy $t_0 <t'<t_i$ and
$a_i<f'(t_i)<b_i$ for all $i=1,\ldots,n$.

\blem
\label{sigid3}
$\f_\h\subset\mathfrak C$
\elem

\noindent{\bf Proof }
It is enough to generate the $\epsilon$-balls in $\cal H$
with cylinders.
We will start with $\epsilon$-balls around points of $\cal H$
consisting of finitely many paths of $\Pi$.

We will use the concept of a {\em cone} in $\R^2$
around $(f,t)$.
Let $r^-=r^-(t,\epsilon)$ and $r^+=r^+(t,\epsilon)$
be the two solutions of
\begin{equation}
{|\tanh(r)-\tanh(t)|}=\epsilon,
\end{equation}
with $r^-\leq r^+$.
For $s$ fixed, let
$x^-(s)=x^-(s,\epsilon)$ and $x^+(s)=x^+(s,\epsilon)$
be the solutions for small $\epsilon$ of
\begin{equation}
\frac{\tanh(x)-\tanh(\hat f(s))}{|s|+1}=\pm\epsilon,
\end{equation}
with $x^-(s)\leq x^+(s)$.
The cone around $(f,t)$ is defined as
\begin{equation}
\C:=\{(x,y)\in\R^2:\,x^-(y)\leq x\leq x^+(y),\,y\geq r^-\}
\end{equation}

Now let $K_0=\{(f_1,t_1),\ldots,(f_n,t_n)\}$. Let
$\C_1,\ldots,\C_n$ be the respective cones of
$(f_1,t_1),\ldots,(f_n,t_n)$. For $i=1,\ldots,n$, let
$r^+_i=r^+(t_i,\epsilon)$, $r^-_i=r^-(t_i,\epsilon)$.

Consider now a family of horizontal lines
$\{L_1,L_2,\ldots\}=\r\times{\cal S}$, where
${\cal S} = \{s_1,s_2, \ldots\}$, with the $s_k$'s distinct and
such that $\cup_{k\geq1}L_k$ is dense in $\R^2$.
For fixed $k$ consider the segments (of nonzero length)
$I_k^i$ into which $L_k$ is divided by all the points of the form
$x_i^-(s_k)$ and $x_i^+(s_k)$
for $i = 1, \ldots ,n$ (number
of such segments $\leq 2n+1$). Segments with interior points in some
cone are closed; otherwise, they are open.

Let ${\cal I} = \{I^{i_1}_{k_1}, \ldots, I^{i_m}_{k_m}\}$ be any
finite sequence of the intervals defined above, with
$k_1,\ldots,k_m$ distinct. For $1 \leq i \leq n$, we will say that
${\cal I}$ is {\em i-good} if $\C_i$ contains all the intervals in
${\cal I}$. If ${\cal I}$ is not $i$-good for any $1 \leq i \leq
n$, then ${\cal I}$ is {\em bad}. Let  \beqn \nonumber \hat
C_i&:=&\{K\in{\cal H}:
\mbox{ there exists $(f,t)\in K$ such that}\\
\label{eq:cil3} &&\quad\quad\quad t\in[r^-_i,r^+_i] \ \mbox{and
$\{f(s)\}\times\{s\} \in \C_i \hbox{ for all } s \geq t$}\} \\
\nonumber &=&
\{K\in \h:  \mbox{ there exists $(f,t)\in K$ such that}\\
\label{eq:cil4} &&\quad\quad\quad d((f,t),(f_i,t_i)) \leq
\epsilon\}. \eeqn
It is not hard to see that $\hat C_i$  belongs to $\mathfrak C$ by
writing $[r_i^-,r_i^+]$ as a finite union of small subintervals 
and then approximating $\hat C_i$ by a finite union of sets of the
form
$\bar C_{{\cal
I}}^{s}$, where: ${\cal I} = \{I_{1}, \ldots, I_{m}\};$  $I_{\ell}
= [x_i^-(s'_\ell),x_i^+(s'_\ell)]\times\{s'_\ell\};$
$s \leq s'_1 \leq s'_2\leq\ldots$; $ s,s'_1 \in [r_i^-,r_i^+];$ and each
 $s'_\ell \in {\cal S}$. Note that in the definition (\ref{eq:cil1})
of such a $\bar C_{{\cal I}}^{s}$, the starting time $t$ of the path
$(f,t)$ must be in $[s,s'_1]$. 
Define $\hat{C}:= \cap_{i=1}^n \hat{C}_i$, $i=1,\ldots,n$. 
Then,
\begin{equation}
\label{eq:rep} \overline{B_{\cal H}(K_0,\epsilon)}=
\left[\left(\bigcup_{m, {\cal I}:\, {\cal I} \mbox{ \tiny is bad}}
C_{{\cal I}}\right) \bigcup
\left(\bigcup_{i=1}^n\bigcup_{\stack{m,{\cal I}, k:\,\, {\cal I}
\mbox{ \tiny is $i$-good}\mbox{ \tiny and}} {s_k>\,\sup_j(r^+_j:\,
{\cal I}\mbox{ \tiny is $j$-good})}} C_{{\cal
I}}^{s_k}\right)\right]^c \bigcap {\hat{C}}.
\end{equation}

We point out that the expression inside the square brackets on the right
hand side is the set of all points $K \in \h$ such that for some
$i=1,\ldots,n$, there is a  path in $K$
 at a distance greater than $\epsilon$ from
$(f_i,t_i)$.  Bad ${\cal I}$'s, in the first term, ensure that
some path in $K$ is at distance greater than $\epsilon$ of
$(f_i,t_i)$ for some $i=1,\ldots,n$ spatially.
 In the second term, some path in $K$
starts at a distance greater than $\epsilon$ from the starting
time of some $(f_i,t_i)$. $\hat{C}$ is the set of all $K\in \h$
with the property, that for all $1 \leq i \leq n$, $K$ contains a
path which is within $\epsilon$ of $(f_i,t_i)$.

To generate  $B_{\cal H}(K,\epsilon)$ for arbitrary
$K\in\cal H$, we approximate $B_{\cal H}(K,\epsilon)$
by an increasing sequence of balls around $\tilde K$'s consisting
of finitely many paths. For that, we note that, by compactness
of $K$, for every integer $j>1$, there exists $K_j\in\cal H$ consisting
of finitely many paths such that $K_j\subset K$ (as subsets
of $\Pi$) and $d_{\cal H}(K,K_j)<\epsilon/j$ for all $j>1$.
We then have
\begin{equation}
B_{\cal H}(K,(1-2/j)\epsilon)\subset B_{\cal
H}(K_j,(1-1/j)\epsilon)\subset
B_{\cal H}(K,\epsilon).
\end{equation}
The first inclusion is justified as follows.
Let $K'\in B_{\cal H}(K,(1-1/j)\epsilon)$. Then
$d_{\cal H}(K,K')<(1-1/j)\epsilon$ and, by the triangle inequality,
\begin{equation}
d_{\cal H}(K_j,K')\leq d_{\cal H}(K_j,K)+d_{\cal H}(K,K')<
\epsilon/j+(1-2/j)\epsilon=(1-1/j)\epsilon.
\end{equation}
Thus
$K'\in B_{\cal H}(K_j,\epsilon)$.
The second inclusion is justified similarly.
It is clear now that
$\cup_{j>1}B_{\cal H}(K,(1-2/j)\epsilon)=
\cup_{j>1}B_{\cal H}(K_j,(1-1/j)\epsilon)=B_{\cal H}(K,\epsilon)$.

%%%%%%%%%%%%%%%%%%%%%%%%%%%%%%%%%%%%%%%%%%%%%%
%%%%%%%%%%%%%% APPENDIX B %%%%%%%%%%%%%%%%%%%%%%%%%
%%%%%%%%%%%%%%%%%%%%%%%%%%%%%%%%%%%%%%%%%%%%%%
%%%%%%%%%%%%%%%%%%%%%%%%%%%%%%%%%%%%%%%%%%%%%%

\section{Compactness and tightness}
\label{app:comptight}
\setcounter{equation}{0}

\begin{figure}[!ht]
\begin{center}
\includegraphics[width=10cm]{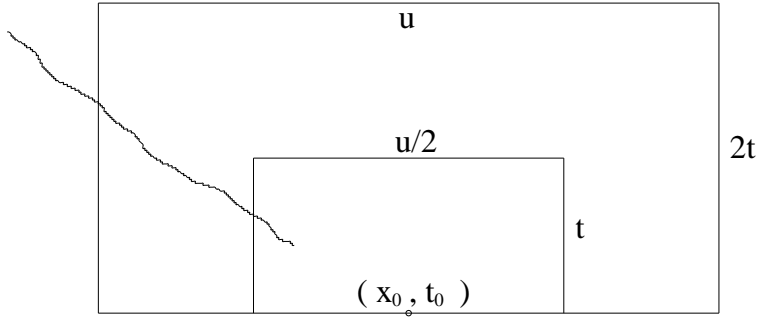}
\caption{Schematic diagram of a path causing the unlikely event
$A_{t,u}(x_0,t_0)$ to occur.} \label{figtight}
\label{tight} %don't know what this does
\end{center}
\end{figure}

Let $\Lambda_{L,T}
= [-L,L]\times [-T,T]$, and let $\{\mu_m\}$ be a sequence of probability
measures
on $(\h,\f_\h)$. For $x_0,t_0\in\r$ and $u,t>0$,
let $R(x_0,t_0;u,t)$ denote the rectangle
$[x_0-u/2,x_0+u/2]\times[t_0,t_0+t]$ in $\R^2$.
Define $A_{t,u}(x_0,t_0)$ to be the event
(in $\f_\h$) that $K$ (in $\h$) contains a path touching both
$R(x_0,t_0;\frac{u}{2},t)$ and (at a later time) the left or right
boundary of the bigger rectangle $R(x_0,t_0;u,2t)$. See
Figure~\ref{figtight}.

Our tightness condition is:
\begin{equation}
\nn
(T_1) \quad \tilde g(t,u;L,T)\equiv
t^{-1}\limsup_m\,\,\sup_{(x_0,t_0) \in \Lambda_{L,T}}
\mu_m(A_{t,u}(x_0,t_0))\to0\mbox{ as }t\to0\!+
\end{equation}
\bprop
\label{prop:tight}
Condition $T_1$ implies tightness of $\{\mu_m\}$.
\eprop
{\bf Proof } Let $$g_m (t,u;L,T) =\sup_{(x_0,t_0)\in\Lambda_{L,T} }
\mu_m(A_{t,u}(x_0,t_0)).$$
Now define $B_{t,u}(x_0,t_0)$ as the event (in $\f_\h$) that
$K$ (in $\h$) contains a path which touches a point $(x',t') =
(f(t'), t') \in
R(x_0,t_0;u/2,t)$ and for some $t'' \in [t',t'+t], |f(t'')-f(t')|
\geq u$. We observe that $B_{t,u}(x_0,t_0)) \subseteq
A_{t,u}(x_0,t_0)$.

We now cover $\Lambda_{L,T}$  with  $\frac{u}{2} \times t$ rectangular
boxes.
Let $L_D =L_D(u)= \{ -L + k\frac{u}{2} : k \in \Z, 0 \leq k \leq
\lceil\frac{2L}{u/2}\rceil\}$
and $T_D = T_D(t)=\{ -T + mt : m \in \Z, 0 \leq m \leq \lceil\frac{2T}{
t}\rceil\}$. Then,
\beqn
&&\mu_m(\cup_{(x_0,t_0)\in \Lambda_{L,T}}B_{t,u}(x_0,t_0))\\
  &\leq & \mu_m(\cup_{(x_0,t_0)\in L_D\times T_D}B_{t,u}(x_0,t_0))
\label{two}
\\
&\leq &  \mu_m(\cup_{(x_0,t_0)\in L_D\times T_D}A_{t,u}(x_0,t_0))\\
&\leq & \left\lceil \frac{2L+1}{u/2}\right\rceil
\left\lceil \frac{2T+1}{t}\right\rceil
g_m(t,u;L,T) \\
& \leq & C' \frac{LT}{tu} g_m(t,u;L,T) \\
& \leq & C' \frac{LT}{u} (\tilde g(t,u;L,T)+\delta).
\label{six}
\eeqn
for any $\delta>0$, where in~(\ref{six}) $m$ is larger than some
$M(t,u;L,T;\delta)$.
The first inequality follows from the observation that
if $K$ (in $\h$) is an outcome in $B_{t,u}(x,t)$ for some $(x,t) \in
   \Lambda_{L,T}$ then $K$ is an outcome in
$B_{t,u}(x',t')$ for some $(x',t') \in L_D\times T_D$.

Now let $\{\phi_n\}$ be a sequence of positive real numbers
with $\lim_{n \to \infty} \phi_n  = 0$.
Now choose $L_n \to \infty$ and $T_n \to \infty$
such that $|\Phi(x,t)| < \frac{\phi_n}{3}$
if $|x| \geq L_n$ or $ t \geq T_n$. Let $\{u_n\}$ be the sequence of
real numbers where $u_n = \frac{\phi_n}{3}$. Now choose sequences of
positive
real numbers $\{t'_n\},\{\delta_n\} \to 0$ such that $C'\frac{L_n T_n}{u_n}
    (\tilde g(t'_n,u_n,L_n,T_n) +\delta_n)\leq
\frac{1}{2^n}$.
For all $n \in \N$, let
\[
C_n(t) = C(t,u_n;L_n,T_n) \equiv \cup_{(x_0,t_0)\in \Lambda_{L_n,T_n}}
B_{t,u_n}(x_0,t_0).
\]
Then, from (\ref{two}) and (\ref{six}) we have that
if $m \geq  M(t'_n,u_n;L_n,T_n;\delta_n)$, then
$\mu_m(C_n(t'_n)) \leq \frac{1}{2^n}$. Since $(\h,d_{\h})$ is a complete
separable metric space, any single measure on $(\h,\f_\h)$ is tight.
Therefore
there exists $D^\epsilon_M $, a compact set in $\h$, such that
$\mu_m(D^\epsilon_M) \geq 1- \frac{\epsilon}{2}$, for all $m \leq M$.

Let $K \in C_n^c$ be a compact set of paths.
   Let $\psi_n = \tanh(T_n+t_n) - \tanh(T_n)$. Suppose $\{(f(t), t):
t \geq t_0\} \in K$.
Now, if $t_0 \leq t_1 \leq t_2$ are times such that $|\Psi(t_2) -
\Psi(t_1)|
\leq \psi_n$, then $|\Phi(f(t_1),t_1) - \Phi(f(t_2),t_2)| \leq \phi_n$.
Let $G_n = \cap_{i=n+1}^\infty C_i^c$. Then, for any $m \geq M$,
\begin{equation}
\mu_m (G_n) = 1 - \mu_m (\cup_{i = n+1}^\infty C_i) \geq
1-\sum_{i = n+1}^\infty \frac{1}{2^i} = 1 - \frac{1}{2^n}.
\end{equation}
Let  $D_n = \cup_{K \in G_n} K$. Then $D_n$ is a family
of equicontinuous functions. By the Arzel\`a-Ascoli theorem, $D_n$ is a
compact subset of $\Pi$. Since $G_n$ is a collection of
   compact
subsets of $D_n$, $G_n$ is a compact subset of $\h$. Let $\epsilon
> 0$. Choose $n \in \N$ such that $\frac{1}{2^n} < \frac{\epsilon}{2}$.
Let $K_\epsilon = D^\epsilon_M \cap G_n$.
Then we have
\[
\sup_m \mu_m (K_\epsilon) \geq 1 - \epsilon
\]
\noindent where  $K_\epsilon$ is a compact subset of $\h$. This proves
that the family of measures $\{\mu_m\}$ is tight.

\brm
\label{holder}
An argument similar to that for Proposition~\ref{prop:tight} can be made to show
that, if instead of $(T_1)$, one has the condition
\begin{itemize}
\item[$(T_1')$]\quad $\sum_{t: t=2^{-k}, k \in\N} t^{-(1+\alpha)}
\sup_m \sup_{x_0,t_0}\mu_m\left(A_{t,t^\alpha}(x_0,t_0)\right) < \infty$
\end{itemize}
for some $\alpha>0$, then each $\mu_m$ as well as any subsequential limit
$\mu$ of $(\mu_m)$ is supported on paths
which are H\"{o}lder continuous with index $\alpha$.
\erm

\bprop
\label{prop: nctightness}
Suppose $\{\X_m\}$ is a sequence of $(\h,d_{\h})$-valued random
variables whose paths are non-crossing. Suppose in addition,

$(I'_1)$ For each $y \in {\cal D}$, there exist (measurable)
path-valued random variables $\theta_m^y \in \X_m$ such that $\theta_m^y$
 converges in
distribution to a Brownian motion $Z_{y}$ starting
at $y$.

\noindent Then the distributions $\{\mu_m\}$ of $\{\X_m\}$ are tight.
\eprop

\noindent {\bf Proof }
 From the proof of the Proposition~\ref{prop:tight}, it is
sufficient to show that for each $u >0$,
\[
\limsup_{m}\mu_m(\cup_{(x_0,t_0)\in L_D(u)\times T_D(t)}B_{t,u}(x_0,t_0))
\to 0
\hbox{ as } t \to 0.
\]

For $u>0, t>0, (x_0,t_0) \in \R^2$, choose two points
 $y_1,y_2 \in {\cal D}$ from the two rectangles $R(x_0\mp \frac{3}{8} u,
t_0-\frac{t}{2};\frac{u}{8},\frac{t}{4})$  respectively. Let
\[
B_1^m(x_0,t_0,t,u) = \left\{ K \in \h \left| \max_{ s \le t_0 + 2t }
|\theta^{y_1}_m(s) -y_1| < u/16\right.\right\}
\]
\[
B_2^m(x_0,t_0,t,u) = \left\{ K \in \h \left| \max_{ s \le t_0 + 2t }
|\theta^{y_2}_m(s) - y_2| < u/16\right.\right\}
\]
and
\(
D^m_{t,u}(x_0,t_0) = B_1^m \cap B_2^m.
\)
Now observe that $D^m_{t,u}(x_0,t_0) \subseteq
B^c_{t,u}(x_0,t_0)$ for large enough $m$. Therefore we have
\beqn
& &\limsup_m \mu_m(\cup_{(x_0,t_0)\in L_D\times T_D}B_{t,u}(x_0,t_0))
\\
&\leq  & \sum_{(x_0,t_0)\in L_D\times T_D} [1 - \liminf_m
\mu_m(D^m_{t,u}(x_0,t_0))]
\eeqn
\noindent Since
$\theta_m^y$
converges in distribution to a Brownian motion
$Z_{y}$ starting at $y$, we have
\beqn
 \liminf_m (\mu_m (B_1^m)) &=&
\P\left(\max_{s \le t_0 + 2t }
|X_{y_1}(s) - y_1| < u/16\right)\\
 & \geq &1 - C t^2/u^4
\eeqn
\noindent and
\beqn
 \liminf_m (\mu_m (B_2^m)) &=&
 \P\left(\max_{ s \le t_0 + 2t }
\left|X_{y_2}(s) - y_2| < u/16\right.\right)\\
&  \geq &1 - C t^2/u^4
\eeqn
\noindent Therefore we have
\[
\liminf_m \mu_m(D^m_{t,u}(x_0,t_0)) \geq 1- 2 C t^2/u^4
\]
\noindent which gives us
\[
\limsup_m \mu_m(\cup_{(x_0,t_0)\in L_D(u)\times T_D(t)}B_{t,u}(x_0,t_0))
\leq 2 C \sum_{(x_0,t_0)\in L_D(u)\times T_D(t)} t^2/u^4
\]
\noindent Since $|L_D(u)\times T_D(t)|\sim \frac{1}{ut}$ we have shown that
\[
\limsup_{m}\mu_m(\cup_{(x_0,t_0)\in L_D(u)\times T_D(t)}B_{t,u}(x_0,t_0))
\to 0
\hbox{ as } t \to 0,
\]
and the proof is complete.

\brm
The proof of Proposition \ref{prop: nctightness} shows that the
limiting processes $Z_y$ starting at $y = (\bar{x},\bar{t})$
need not be Brownian motions.
It is sufficient that they be continuous processes such that
for each fixed $u>0$,
\begin{equation}
\label{gb1}
\frac{1}{t} \sup_{y} \P(\sup_{ \bar{t} \leq s \leq \bar{t} + t} |
Z_y(s) - Z_y(\bar{t})| \geq u) \to 0 \hbox{ as } t \to 0^+.
\end{equation}
\erm

\bprop
\label{prop:compact1}
Let ${\cal D}$ be a countable dense subset of $\R^2$ and let
 $\mu_k$ be the distribution of the $(\h,\f_{\h})$-valued random
variable ${\cal W}_k ={\cal W}_k({\cal D}) = \{ \tilde{{\cal W}_1}, \ldots ,
\tilde{{\cal W}_k} \}$.
Then
the family of measures $\{\mu_k\}$ is tight.
\eprop
\noindent {\bf Proof } This is an immediate consequence of
Proposition(\ref{prop: nctightness}).

\bprop
\label{prop:tight-compact}
If $\W_n$ is an a.s. increasing sequence of $(\h,d_\h)$- valued
random variables and the family of  distributions $\{\mu_n\}$ of $\W_n$ is tight,
then $\overline{\cup_n \W_n} $is almost surely compact
 (in $(\Pi,d)$).
\eprop

{\bf Proof } Let $\tilde{\W_k}$ be an increasing sequence of points
(subsets of $\Pi$) in $(\h,d_\h)$, which converge in $d_\h$ metric
 to some
point $\tilde{\W}$ in $(\h,d_\h)$. If
for some $k$, $\tilde{\W_k}$ is not a subset of $\tilde{\W}$, then
there exists an $\e>0$ such that $d_\h(\tilde{\W},\tilde{\W_n}) >\e$
for all $n \geq k$ contradicting the claim that $\tilde{\W_k}$ converges
to $\tilde{\W}$. Therefore $\tilde{\W_k} \subseteq \tilde{\W}$ for all $k$.
This implies $ \overline{\cup_k \tilde{\W_k}}\subseteq\tilde{\W}$ and
therefore is a compact subset of $\Pi$ since it is a closed subset
of the compact set $\tilde{\W}$. Since
$\{\mu_n\}$ is tight, given an $\e > 0$, there exists a compact
subset K of $\h$ such that $\P(\W_n \in K) \geq
1 - \e$ for all $n$, so by monotonicity,
$P(\W_n \in K$ for all $n)\, \geq \, 1 - \e$.
But if $\W_n \in K$ for all $n$, then since $K$
is compact,
there exists a subsequence $\W_{n_j}$ which converges to a point in
$K$ and thus in $\h$. This implies by the first part of this proof
that $\overline{\cup_{n_j} \W_{n_j}}$($ =
\overline{\cup_n \W_n}$ because $\W_n$ is increasing in $n$)
 is a compact subset of $\Pi$. Thus we
have shown that $\P(\overline{\cup_n \W_n}
\hbox{ is a compact subset of } \Pi) \geq 1-\e$. Since the claim
is true for all $\e >0$ we have proved the proposition.

\bprop
\label{prop: compact2}
Let $\hat{{\cal D}} = \{ (\hat{x_i},\hat{t_i}): i = 1,2, \ldots\}$
be a (deterministic) dense countable subset of $\R^2$ and let $\{
\hat{{\cal W}_i}:i= 1,2, \ldots\}$ be $(\Pi,d)$-valued random
variables starting from $(\hat{x_i},\hat{t_i})$. Suppose that
the joint distribution of each finite subset of
the $\hat{\cal W}_i$'s is that of
coalescing Brownian motions. Then $\overline{\cup_{n=1}^\infty
\{\hat{{\cal W}_1}, \hat{{\cal W}_2} \ldots   \hat{{\cal W}_n}\}}$
is almost surely compact. In particular for ${\cal W}_n$ defined in
Proposition (\ref{prop:compact1}),
$\bw = \overline{\cup_n {\cal W}_n} $is almost surely compact.
\eprop

{\bf Proof } The proof follows immediately from Propositions
\ref{prop:compact1} and \ref{prop:tight-compact}.

\bprop
\label{prop: D_1-D_2}
Let $\hat{{\cal D}}$ and $\{\hat{{\cal W}_i}:i= 1,2, \ldots\}$
be as in Proposition~\ref{prop: compact2}
and let $\{ \hat{{\cal W'}_i}:i= 1,2, \ldots\}$ (on some other probability
space) be equidistributed with $\{\hat{{\cal W}_i}:i= 1,2, \ldots\}$. Then
$\hat{{\cal W}} \equiv \overline{\{{\hat{\cal W}_i}:i= 1,2, \ldots\}}$ and
$\hat{{\cal W'}} \equiv \overline{\{ {\hat{\cal W'}_i}:i= 1,2, \ldots\}}$
are equidistributed $(\h,\f_\h)$-valued random variables.
\eprop
{\bf Proof } It is an easy consequence of \ref{prop: compact2} that
$\{\hat{{\cal W}_i}:i= 1, \ldots,n\}$ (respectively
$\{ \hat{{\cal W'}_i}:i= 1, \ldots,n\}$) converges a.s.
as $n \to \infty$ in $(\h,d_\h)$ to
$\hat{{\cal W}}$ (respectively $\hat{{\cal W'}}$). But then the identical
distributions of $\{\hat{{\cal W}_i}:i= 1, \ldots,n\}$ and
$\{ \hat{{\cal W'}_i}:i= 1, \ldots,n \}$ converge respectively to the
distributions of $\hat{{\cal W}}$ and  $\hat{{\cal W'}}$, which thus
must be identical.

%%%%%%%%%%%%%%%%%%%%%%%%%%%%%%%%%%%%%%%%%%%%%%%
%%%%%%%%%%%%%%% APPENDIX C %%%%%%%%%%%%%%%%%%%%%%%%%
%%%%%%%%%%%%%%%%%%%%%%%%%%%%%%%%%%%%%%%%%%%%%%%
%%%%%%%%%%%%%%%%%%%%%%%%%%%%%%%%%%%%%%%%%%%%%%%

\section{Hausdorff dimension of the graph of the sum of two Brownian motions}
\label{app:haus}
\setcounter{equation}{0}

\bprop
\label{prop:haus}
Let $X,Y$ be two independent standard Brownian motions and let ${\cal T}^+$ denote the set
of record times of $X$, i.e., ${\cal T}^+=\{t\geq0:\,X(t)=M(t)\}$,
where $M(t):=\sup_{0\leq s\leq t}X(s)$ is the maximum of X up to time $t$.
Then, for $t_0>0$ and $a,b\in\r$ with $|a|+|b|>0$, the set
${\cal G}^+:=\{(aX(t)+bY(t),t):\,t\in{\cal T}^+\cap[0,t_0]\}$,
the projection of ${\cal T}^+\cap[0,t_0]$ onto the graph of $aX+bY$,
has Hausdorff dimension $1$ almost surely.
\eprop

\noindent{\bf Proof } An upper bound of $1$ for the
Hausdorff dimension follows readily from the
fact that ${\cal G}^+$ is the image of a set, ${\cal T}^+\cap[0,t_0]$,
of Hausdorff dimension $1/2$ a.s.~(since ${\cal T}^+\cap[0,t_0]$ has the same
distribution as the set of zeros of $X$, $\{t\in[0,t_0]:\,X(t)=0\}$;
this follows from $(M(t)-X(t):\,0\leq t\leq t_0)$ having the same distribution
as $(|X(t)|:\,0\leq t\leq t_0)$~\cite{kn:L}) --- see~\cite{kn:Ta} --- under a
map which is a.s.~(uniformly) H\"older continuous of exponent $\alpha$ for every
$\alpha<1/2$, namely, the map $t\to(aX(t)+bY(t),t)$, where we use the well
known H\"older continuity properties of Brownian motion.

The desired lower bound is obtained by noting that the Hausdorff dimension of
${\cal G}^+$ is bounded below by the Hausdorff dimension of the image of
${\cal T}^+\cap[0,t_0]$ under $aX+bY$, or equivalently under $aM+bY$, namely
$\{aM(t)+bY(t):\,t\in{\cal T}^+\cap[0,t_0]\}$. Notice that the latter set
equals
$\{as+bY(T(s)):\,s\in[0,M(t_0)]\}$,
where $T$ is the hitting time process associated to $X$, defined as
$T(x) := \inf\{t\geq0:\,X(t)=x\}$.
It suffices to show that the Hausdorff dimension of
$\{as+bY(T(s)):\,s\in[0,L]\}$
is a.s.~greater than or equal to $1$ for every deterministic $L>0$.
But that follows from known results as well. $Z(t):=at+bY(T(t))$ is a self similar process of
exponent $1$ with stationary increments and satisfies also the following condition
of Theorem 3.3 in~\cite{kn:XL}, from which the dimension bound follows. The condition is that
there exists a constant $K$ such that $\P(|Z(1)|\leq x)\leq Kx$ for every $x\geq0$.
This property is readily obtained from the distributions of
$Y$ and the hitting time variable $T(1)$.

%%%%%%%%%%%%%%%%% REFERENCES %%%%%%%%%%%%%%%%%%%

\noindent L.~R.~G.~Fontes \\
Instituto de Matem\'atica e Estat\'{\i}stica, Universidade de S\~ao
Paulo\\ Rua do Mat\~ao 1010, CEP 05508--090, S\~ao Paulo SP,
Brazil

\vspace{.5cm}

\noindent M.~Isopi\\
Dipartimento di Matematica, Universit\`a di Bari\\
Via E.~Orabona 4, 70125 Bari, Italia

\vspace{.5cm}

\noindent C.~M.~Newman\\ 
Courant Institute of Mathematical Sciences, New York University\\
New York, NY 10012, USA

\vspace{.5cm}

\noindent K.~Ravishankar\\
Department of Mathematics, SUNY-College at New Paltz \\
New Paltz, New York 12561, USA


\begin{thebibliography}{xxxxxx 89}


\bibitem{kn:FINR}  L.~R.~G.~Fontes, M.~Isopi, C.~M.~Newman, K.~Ravishankar,
                   The Brownian web,
                   {\em Proc.~Nat.~Acad.~Sciences} {\bf 99}, 15888-15893 (2002).

\bibitem{kn:A2}   R.~Arratia,
                  Coalescing Brownian motions and the voter model on $\Z$,
                  Unpublished partial manuscript (circa 1981),
                  available from rarratia@math.usc.edu.

\bibitem{kn:TW}   B.~T\' oth, W.~Werner,
                  The true self-repelling motion,
                  {\em Probab. Theory Related Fields} {\bf 111}, 375-452 (1998).

\bibitem{kn:FINS1} L.~R.~G.~Fontes, M.~Isopi, C.~M.~Newman, D.L.~Stein,
                  Aging in 1${}{D}$ discrete spin models and equivalent systems,
                  {\em Phys. Rev. Lett.} {\bf 87}, 110201-1 -- 110201-4 (2001).

\bibitem{kn:FINS2} L.~R.~G.~Fontes, M.~Isopi, C.~M.~Newman, D.L.~Stein,
                  1${}{D}$ Aging,
                  in preparation (2003).

\bibitem{kn:D}    M.~D.~Donsker,
                  An invariance principle for certain probability limit theorems,
                  {\em Memoirs of the AMS} {\bf 6}, 1-12 (1951).

\bibitem{kn:H}    T.~E.~Harris,
                  Additive set-valued Markov processes and graphical methods,
                  {\em Ann.~Probability} {\bf 6}, 355-378 (1978).

\bibitem{kn:FFW}  P.~A.~Ferrari, L.~R.~G.~Fontes, X.~Y.~Wu,
                  preprint.

\bibitem{kn:S}    A.~E.~Scheidegger,
                  A stochastic model for drainage patterns into an intramontane trench,
                  {\em Bull.~Ass.~Sci.~Hydrol.}~{\bf 12}, 15-20 (1967).

\bibitem{kn:RR}   I.~Rodriguez-Iturbe, A.~Rinaldo,
                  {\em Fractal river basins: chance and self-organization},
                  Cambridge Univ.~Press, New York (1997).

\bibitem{kn:A0}   R.~Arratia,
                  {\em Coalescing Brownian motions on the line,}
                  Ph.D. Thesis, University of Wisconsin, Madison (1979).

\bibitem{kn:A1}   R.~Arratia,
                  Limiting point processes for rescalings of coalescing and
                  annihilating random walks on $\Z\sp{d}$,
                  {\em Ann.~Prob.}~{\bf 9}, 909-936 (1981).

\bibitem{kn:H2}   T.~E.~Harris,
                  Coalescing and noncoalescing flows in $\R^1$,
                  {\em Stoch. Proc. and their Appl.}~{\bf 17}, 187-210 (1984).

\bibitem{kn:STW}  F.~Soucaliuc, B.~T\' oth, W.~Werner,
                  Reflection and coalescence between independent one-dimensional Brownian paths,
                  {\em Ann.~Inst.~H.~Poincar\' e Probab.~Statist.}~{\bf 36}, 509-545 (2000).

\bibitem{kn:T}    B.~Tsirelson,
                  White noises, black noises and their scaling limits,
                  2002 unpublished lecture notes, available from\\
                  arXiv: math.PR/0301237

\bibitem{kn:A}    M.~Aizenman,
                  Scaling limit for the incipient spanning clusters,
                  pp.~1-24 in {\em Mathematics of Multiscale Materials:
                  Percolation and Composites (Minneapolis, Minn., 1995-1996)},
                  IMA Vol.~Math.~Appl.~{\bf 99}, Springer-Verlag, New York (1998).

\bibitem{kn:AB}   M.~Aizenman, A.~Burchard,
                  H\"older regularity and dimension bounds for random curves,
                  {\em Duke Math.~J.}~{\bf 99}, 419-453 (1999).

\bibitem{kn:ABNW} M.~Aizenman, A.~Burchard, C.~M.~Newman, D.~B.~Wilson,
                  Scaling limits for minimal and random spanning trees in two dimensions,
                  {\em Random Structures Algorithms}~{\bf 15}, 319-367 (1999).

\bibitem{kn:V}    S.~R.~S.~Varadhan,
                  {\em Stochastic Processes},
                  Courant Inst.~of Math.~Sciences, New York (1968).

\bibitem{kn:P}    V.~V.~Piterbarg,
                  Expansions and contractions of isotropic stochastic flows of homeomorphisms,
                  {\em Ann.~Probability} {\bf 26}, 479-499 (1998).

\bibitem{kn:BG}   M.~Bramson, D.~Griffeath,
                  Clustering and dispersion rates for some interacting
                  particle systems on $\Z^1$,
                  {\em Ann.~Probability} {\bf 8}, 183-213 (1980).

\bibitem{kn:L}    P.~L\'evy,
                  {\em Processus Stochastiques et Mouvement Brownien},\\
                  Gauthier-Villars, Paris (1948).

\bibitem{kn:Ta}   S.~J.~Taylor,
                  The $\alpha$-dimensional measure of the graph and set of zeros of a Brownian path,
                  {\em Proc.~Cambridge Philos.~Soc.}~{\bf 51}, 265--274 (1955).

\bibitem{kn:XL}   Xiao~Yimin, Lin~Huonan,
                  Dimension properties of sample paths of self-similar processes,
                  {\em Acta Math.~Sinica (N.S.)}~{\bf 10}, 289--300 (1994).

\end{thebibliography}
\end{document}